\documentclass{article}
\usepackage{arxiv}

\usepackage{amsmath,amsfonts,amssymb}
\usepackage{mathtools}
\usepackage[hidelinks]{hyperref}
\usepackage{algorithmic}
\usepackage{algorithm}
\usepackage{array}
\usepackage[caption=false,font=normalsize,labelfont=sf,textfont=sf]{subfig}
\usepackage{relsize}
\usepackage{textcomp}
\usepackage{stfloats}
\usepackage{url}
\usepackage{verbatim}
\usepackage{graphicx}
\usepackage{natbib}
\usepackage{doi}
\usepackage{xcolor}
\usepackage{paralist}
\usepackage{blindtext}
\usepackage{multicol}
\usepackage{multirow}
\usepackage[subnum]{cases}

\newtheorem{theorem}{Theorem}[section]

\newtheorem{remark}[theorem]{Remark}

\newcount\Comments
\usepackage{color}
\newcommand{\kibitz}[2]{\ifnum\Comments=0\textcolor{#1}{#2}\fi}

\title{Cyberattacks on Adaptive Cruise Control Vehicles: An Analytical Characterization}

\author{Shian Wang\\
	Electrical and Computer Engineering\\
	The University of Texas at El Paso\\
	\texttt{swang14@utep.edu}\\
	\And
	Mingfeng Shang\\
	Civil, Environmental, and Geo-Engineering\\
	University of Minnesota\\
	\texttt{shang140@umn.edu} \\
	\And
	Raphael Stern\\
	Civil, Environmental, and Geo-Engineering\\
	University of Minnesota\\
	\texttt{rstern@umn.edu} \\
}

\renewcommand{\shorttitle}{Analytical Characterization of Cyberattacks on ACC Vehicles}

\begin{document}\sloppy
\maketitle

\begin{abstract}
    While automated vehicles (AVs) are expected to revolutionize future transportation systems, emerging AV technologies open a door for malicious actors to compromise intelligent vehicles. As the first generation of AVs, adaptive cruise control (ACC) vehicles are vulnerable to cyberattacks. While recent effort has been made to understanding the impact of attacks on transportation systems, little work has been done to systematically model and characterize the malicious nature of candidate attacks. In this study, we develop a general framework for modeling and synthesizing two types of candidate attacks on ACC vehicles, namely direct attacks on vehicle control commands and false data injection attacks on sensor measurement, with explicit characterization of their adverse effects. Based on linear stability analysis of car-following dynamics, we derive a series of analytical conditions characterizing the malicious nature of potential attacks. This ensures a higher degree of realism in modeling attacks with adverse effects, as opposed to simply considering attacks as constants or random variables. Notably, the conditions derived provide an effective method for strategically synthesizing an array of candidate attacks on ACC vehicles. We conduct extensive simulation to examine the impacts of intelligently designed attacks on microscopic car-following dynamics and macroscopic traffic flow. Numerical results illustrate the mechanism of candidate attacks, offering useful insights into understanding the vulnerability of future transportation systems. The methodology developed allows for further study of the widespread impact of strategically designed attacks on traffic cybersecurity, and is expected to inspire the development of efficient attack detection techniques and advanced vehicle controls. 
\end{abstract}


\section{Introduction}\label{section1}

Automated vehicles (AVs) are expected to revolutionize future transportation systems with considerable benefits, such as reduced energy consumption~\citep{sun2022energy}, improved smoothness of traffic flow~\citep{wang2022optimal}, and optimized parking space allocation~\citep{wang2021optimal}, among others. While emerging AV technologies have a lot to offer in reshaping the future of transportation systems, they open a door for malicious actors to compromise vehicle safety and security~\citep{parkinson2017cyber} which could result in adverse impacts on traffic flow~\citep{dong2020impact}. As the first generation of AVs, adaptive cruise control (ACC) vehicles are intrinsically vulnerable to potential cyberattacks. Those attacks can not only be launched in various forms, but also cause substantial disruption to normal traffic flow~\citep{petit2014potential}. As ACC vehicles are increasingly equipped with advanced vehicle-to-vehicle (V2V) communication technologies, enabling cooperative ACC (CACC) systems, they become even more vulnerable to attacks since V2V communication channels are subject to various adverse actions~\citep{alipour2020impact}. 

Although a cyberattack may only cause subtle changes to the driving behavior of an individual vehicle~\citep{li2018influence}, the resulting impact could propagate along the traffic flow leading to increased traffic congestion and energy consumption of vehicles~\citep{dong2020impact,li2023exploring}. While attacks can be launched in various forms, a few types have drawn particular interest in recent years. For example, cyberattacks can occur to ACC control commands, where the desired acceleration is directly altered by malicious actors~\citep{li2022detecting,zhou2022robust,wang2023minmax}. Another type is false data injection attacks on ACC sensor measurement, which indirectly cause a vehicle to execute undesired maneuvers~\citep{wang2020modeling,wang2023anomaly,wang2023minmax}. In addition, potential adversaries may attempt to shut down a network by temporarily or indefinitely flooding the communication traffic resulting in communication delays~\citep{zeadally2012vehicular}. Such type of attacks, called denial-of-service attacks~\citep{zeadally2012vehicular,wang2020modeling}, are more likely to occur to CACC systems that rely heavily on V2V communication. The reader is referred to~\citep{petit2014potential} for a comprehensive introduction of different forms of potential cyberattacks on intelligent vehicles. It is easy to agree that undesired behavior of attacked ACC vehicles can degrade their own performance. However, it is less apparent how compromised vehicles may fundamentally impact the bulk traffic involving both ACC and human-driven vehicles (HDVs). In fact, numerical simulation has shown that even slight attacks on a single vehicle could lead to reduced traffic throughput and increased energy consumption~\citep{dong2020impact,li2023exploring}, due to the fact that traffic flow results from the collective behavior of all vehicles. Clearly, potential attacks on ACC vehicles pose a significant threat not only to the individual vehicles compromised but also to the safety and efficiency of the entire transportation system. 

While an increasing amount of research effort has been made to understanding the impact of potential attacks on transportation systems in recent years, there is a lack of work that systematically models the malicious nature of attacks with explicit characterization of their adverse effects. Specifically, attacks on ACC vehicles or AVs have been assumed to be constant or stochastic (e.g., normally distributed) in prior studies~\citep{wang2021resilient,wang2020modeling,li2018influence,li2022detecting}, which may not suffice to describe the malicious characteristics of candidate attacks. Instead of launching attacks that are poorly chosen, the attacker is more likely to design candidate attacks in a strategic manner in order to degrade the performance of attacked vehicles thereby compromise the efficiency and safety of the bulk traffic~\citep{wang2023minmax,wang2023novel}. Moreover, synthesizing attacks in a random fashion may result in compromised vehicles being easily detected and identified due to intrusion detection systems~\citep{khraisat2019survey} and anomaly detection techniques~\citep{li2022detecting}. Consequently, questions naturally arise on the realism of modeling cyberattacks on ACC vehicles in the way that has been presented in many relevant studies. In this study, we aim to model and analytically synthesize two types of malicious attacks on ACC vehicles based on the linear stability analysis of car-following dynamics. Notably, we derive a series of mathematical conditions that provide an effective way of synthesizing candidate attacks with adverse effects on traffic flow, in terms of efficiency and safety. Further, we examine the impacts of such strategically designed attacks on microscopic car-following dynamics and macroscopic traffic flow. The main contributions of this study are listed below:
\begin{itemize}
    \item We develop a general framework for modeling and synthesizing two types of candidate attacks on ACC vehicles with explicit mathematical characterization of their malicious nature. The first type is attacks on ACC control commands, while the second one is false data injection attacks on ACC sensor measurement, both of which could result in undesired driving behavior of the attacked vehicles thereby compromising traffic flow. 
    
    \item We derive a series of analytical conditions, based on linear stability analysis of car-following dynamics, for characterizing the malicious nature of potential attacks. This ensures a higher degree of realism in modeling attacks with adverse effects, as opposed to simply considering attacks as constants or random variables. 

    \item We conduct extensive simulation to examine the impacts of strategically designed attacks on microscopic car-following dynamics and macroscopic traffic flow. Numerical results illustrate the mechanism of candidate attacks and offer useful insights into understanding the vulnerability of future transportation systems.

    \item The analysis conducted provides an effective method for synthesizing an array of strategic attacks, thereby allowing for further study of the widespread impact of such strategically designed attacks on future transportation systems. Consequently, it is expected to inspire the development of efficient attack detection techniques and advanced vehicle control algorithms. 
    
\end{itemize}

The remainder of this article is structured as follows. In Section~\ref{section2}, we present a mathematical framework for describing car-following dynamics of mixed traffic, under which two types of candidate attacks are introduced. We provide a series of analytical characterization of the candidate attacks in Sections~\ref{section3} and~\ref{section4}, allowing for attack synthesis in a strategic manner. To illustrate the mechanism of attack synthesis, we conduct extensive numerical simulation in Section~\ref{section5} to examine the impacts of candidate attacks on microscopic car-following dynamics and macroscopic traffic flow. The article is concluded in Section~\ref{section6} with discussion on its limitations and possible extensions.

\section{Mixed Traffic Under Cyberattack}\label{section2}

\subsection{Mathematical model of mixed traffic}\label{section2.1}
\begin{figure}[t!]
	\centering
	\includegraphics[width=\textwidth]{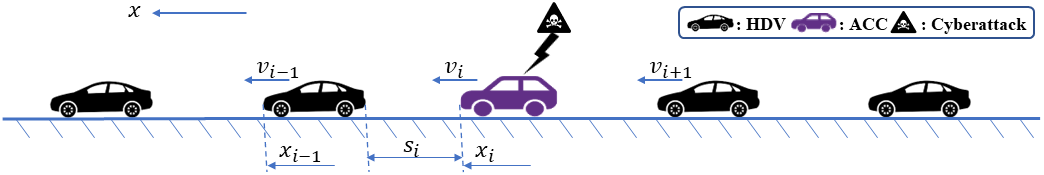}
	\caption{A graphic illustration of mixed traffic involving ACC vehicles and HDVs, with ACC vehicles under potential cyberattacks}
	\label{cyberattack_mixed_traffic}
\end{figure}

We consider mixed traffic involving HDVs and ACC vehicles (or AVs), which is projected to be the case in the foreseeable future~\citep{wang2022optimalTRC}. Following prior studies~\citep{dong2020impact,wang2020modeling,li2022detecting,wang2023anomaly,zhou2022robust} we only consider longitudinal vehicle dynamics, while lateral dynamics could also be studied. Without loss of generality, we consider a string of $m$ vehicles denoted by the totally ordered set $\mathcal{M} = \left\{ 1, 2, 3, \cdots, m \right\}$, $m > 1$, $m \in \mathbb{N}^{+}$. For any vehicle $i \in \mathcal{M}$, let $x_{i}(t)$ and $v_{i}(t)$ denote its position and speed at time $t$, respectively. Consequently, the inter-vehicle spacing between vehicle $i$ and its preceding vehicle $i-1$ is given by $s_{i}(t) = x_{i-1}(t) - x_{i}(t) - l_{i-1}$ with $l_{i-1}$ the length of vehicle $i-1$. These standard notations are commonly used for describing car-following dynamics~\citep{wilson2011car}, with a graphic illustration shown in Fig.~\ref{cyberattack_mixed_traffic}.

Based on the law of physics, the dynamics of any vehicle $i$ is given by the following set of differential equations~\citep{wilson2011car}
\begin{eqnarray}
    &~&\dot{x}_{i}(t) = v_{i}(t),    \label{eq2.1}   \\
    &~&\dot{v}_{i}(t) = f(s_{i}(t), \Delta v_{i}(t), v_{i}(t)), \label{eq2.2}
\end{eqnarray}
where the dot operator denotes differentiation with respect to time; the nonlinear operator $f$ relates acceleration of the $i$th vehicle, $\dot{v}_{i}$, to the variables $s_{i}$, $\Delta v_{i}$ and $v_{i}$; and the relative speed between vehicle $i$ and vehicle $i-1$ is defined as
\begin{equation}
    \Delta v_{i}(t) = \dot{s}_{i}(t) = v_{i-1}(t) - v_{i}(t). \label{eq2.3}
\end{equation}
The above equations~\eqref{eq2.1}--\eqref{eq2.2} represent a generic functional form of car-following dynamics widely adopted in the literature. While the variables involved, including $s_{i}$, $\Delta v_{i}$ and $v_{i}$, are time-variant, for brevity we omit the argument $t$ wherever appropriate. It is worth noting that ACC systems do not require V2V communication. In fact, the measurements, like inter-vehicle spacing $s_{i}$ and relative speed $\Delta v_{i}$, can be readily obtained via ACC onboard radar sensors.

Let $\mathcal{H}$ and $\mathcal{A}$ denote the totally ordered set of HDVs and ACC vehicles in the mixed traffic, respectively. Since human drivers tend to exhibit different driving behaviors from ACC vehicles~\citep{ye2022car}, we further distinguish the functional $f$ of equation~\eqref{eq2.2} as $f_{\text{HDV}}$ and $f_{\text{ACC}}$ for HDVs and ACC vehicles, respectively. The resulting mixed traffic dynamics is given by
\begin{equation}
    \dot{x}_{i} = v_{i}, ~\forall~ i \in \mathcal{M} = \mathcal{H} \cup \mathcal{A}   \label{eq2.4}
\end{equation}
\begin{numcases}{\dot{v}_{i} =}
    f_{\text{HDV}}(s_{i}, \Delta v_{i}, v_{i}), ~\forall~ i \in \mathcal{H}    \label{eq2.5a}  \\  
    f_{\text{ACC}}(s_{i}, \Delta v_{i}, v_{i}), ~\forall~ i \in \mathcal{A}    \label{eq2.5b}
\end{numcases}
It is worth noting that the functionals $f_{\text{HDV}}$ and $f_{\text{ACC}}$ do not need to be same. Instead, they depend on the specific car-following principles abided by HDVs and ACC vehicles. For simplicity, reaction-time delays are not considered in equations~\eqref{eq2.5a}--\eqref{eq2.5b}, as seen in many relevant studies modeling traffic flow based on car-following dynamics~\citep{talebpour2016influence,wang2020modeling,li2022detecting,wilson2011car,wang2022optimal,wang2023general}.

\begin{figure}[t!]
	\centering
	\includegraphics[width=0.5\textwidth]{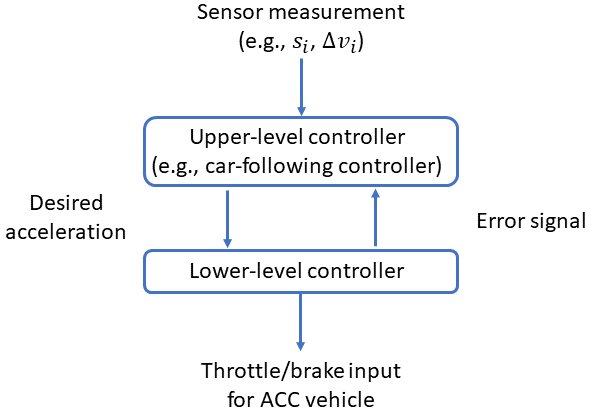}
	\caption{Illustration of two-level controllers for an ACC vehicle, where the upper-level controller determines the desired acceleration while the lower-level controller executes the corresponding throttle or brake input~\citep{rajamani2011vehicle}}
	\label{ACC_controller}
\end{figure}

\subsection{A class of candidate attacks on ACC control commands}\label{section2.2}

The control mechanism of ACC vehicles follows a two-level controller design as shown in Fig.~\ref{ACC_controller}. Specifically, the upper-level controller determines the desired acceleration based on car-following principle like a constant time-gap policy, while the lower-level controller calculates the corresponding throttle or brake input~\citep{rajamani2011vehicle}. In this study, we focus on potential attacks that could occur to the upper-level controller, as seen in prior work~\citep{wang2020modeling,li2022detecting,wang2022planning,wang2023minmax}. We first introduce a class of cyberattacks on ACC control commands, i.e., desired acceleration, termed Type~I attacks, which can pose a significant threat to compromised ACC vehicles, thereby degrading the performance of traffic flow. In prior studies, cyberattacks on intelligent vehicles have been assumed to be deterministic or stochastic with a known probability distribution, without much explicit characterization of the malicious actors~\citep{petit2014potential,li2018influence,wang2020modeling,van2019real,wang2021resilient,li2022detecting}. In contrast, we consider attacks as a bounded unstructured (distribution-free) process that could be deterministic or stochastic, relaxing the assumptions seen in many previous works. 

To this end, we introduce the dynamics of an ACC vehicle under attack. With ACC control commands being attacked, the equivalent effect is the alteration to ACC acceleration from a modeling standpoint~\citep{li2022detecting,wang2022planning,zhou2022robust,wang2023minmax}. Letting $\omega_{i}$ denote the resulting attack on ACC acceleration, the compromised ACC dynamics is then given by
\begin{eqnarray}
    \dot{v}_{i} = f_{\text{ACC}}(s_{i}, \Delta v_{i}, v_{i}) + \omega_{i}, ~\text{if}~ i \in \mathcal{A} ~\text{attacked}.   \label{eq2.6}
\end{eqnarray}
Since intelligent attacks are normally launched based on certain information of the target perceived by the adversaries~\citep{boem2017distributed}, it is reasonable to assume that $\omega_{i}$ is a function of ACC speed. In other words, $v_{i}$ is assumed to be accessed by the malicious actor when launching attacks. Specifically, let $\omega_{i} = g(v_{i})$, where $g: \mathbb{R} \longrightarrow \mathbb{R}$ is an operator denoting the resulting attack on ACC acceleration.

It is noted that equation~\eqref{eq2.6} considering cyberattacks is subject to the physical bounds of vehicle acceleration. Since compromised vehicles do not have prior knowledge of the potential attacks, one can introduce an unstructured process to describe the operator $g$ as follows 
\begin{eqnarray}
    \omega_{i} = g(v_{i}) = \delta_{i}v_{i}, ~\text{if}~ i \in \mathcal{A} ~\text{attacked},   \label{eq2.7}
\end{eqnarray}
where $\delta_{i}$ is an unstructured variable (allowing for inclusion of deterministic and stochastic attack behaviors) assumed to be bounded by some known bounds~\citep{boem2017distributed}, e.g., $\delta_{i} \in \Lambda \coloneqq \left\{\delta_{i}: \left|\delta_{i}\right| \leq r, r \in \mathbb{R}_{\geq 0}\right\}$, with $r=0$ indicating no attack. This covers a broad class of candidate attacks that could potentially occur to ACC control commands. Clearly, a larger value of $r$ indicates more severe attacks that may alter ACC driving behavior to a greater extent. This assumption is reasonable since realistic adversaries generally tend to remain stealthy while being subject to limited resources and energy~\citep{biron2018real}. As a result, the acceleration dynamics of an ACC vehicle under attack is written as
\begin{eqnarray}
    \dot{v}_{i} = f_{\text{ACC}}(s_{i}, \Delta v_{i}, v_{i}) + \delta_{i}v_{i} \coloneqq F_1(s_{i}, \Delta v_{i}, v_{i}, \delta_{i}).   \label{eq2.8}
\end{eqnarray}
Consequently, the mixed traffic dynamics of vehicle acceleration under cyberattack is given by
\begin{numcases}{\dot{v}_{i} =}
    f_{\text{HDV}}(s_{i}, \Delta v_{i}, v_{i}), ~\text{if}~ i \in \mathcal{H},  \label{eq2.9a}  \\
    f_{\text{ACC}}(s_{i}, \Delta v_{i}, v_{i}), ~\text{if}~ i \in \mathcal{A} ~\text{not attacked},    \label{eq2.9b}  \\ 
    F_1(s_{i}, \Delta v_{i}, v_{i}, \delta_{i}), ~\text{if}~ i \in \mathcal{A} ~\text{attacked}.   \label{eq2.9c}
\end{numcases}

\begin{remark}\label{Remark2.1}
Limited prior studies have assumed attacks to be deterministic or stochastic with a given probability distribution like the Gaussian distribution. Unfortunately, this is far from realistic due to the lack of knowledge of malicious actors. In contrast, we characterize candidate attacks with a bounded unstructured process $\left\{\delta_{i}\right\}$, which is capable of describing both deterministic and stochastic processes without being subject to any specific statistical distribution. In addition, we relate attacks to the state (i.e., speed) of an ACC vehicle to capture the intention of adversaries launching attacks based on ACC information perceived. This allows for a higher degree of realism compared to prior studies assuming constant attacks or random attacks using Gaussian noise without explicitly characterizing the intention of malicious actors. Synthesizing attacks with a higher level of realism offers greater insights into understanding their realistic impact on traffic flow, paving the way for developing effective attack detection and mitigation strategies in future studies.
\end{remark}

\subsection{A class of candidate attacks on ACC sensor measurement}\label{section2.3}

While attacks could occur to the desired acceleration as introduced in Section~\ref{section2.2}, they may also corrupt ACC hardware sensors or the software algorithms used for data acquisition in the form of false data injection attacks. It is observed from equation~\eqref{eq2.5b} that the desired acceleration of an ACC vehicle (i.e., the upper-level controller shown in Fig.~\ref{ACC_controller}) depends on sensor measurements like the inter-vehicle spacing $s_i$ and the relative speed $\Delta v_i$. These measurements could be corrupted by malicious actors when used to execute the acceleration command~\citep{wang2020modeling,wang2022planning,wang2023minmax}. Hence, we introduce a class of cyberattacks on ACC sensor measurements, termed Type~II attacks. For any vehicle $i \in \mathcal{A}$, let $\omega_{1,i}$ and $\omega_{2,i}$ denote the false data injection attacks on spacing and relative speed, respectively. Thus, the resulting acceleration dynamics of any attacked ACC vehicle becomes
\begin{eqnarray}
\dot{v}_{i} = f_{\text{ACC}}(s_{i}+\omega_{1,i}, \Delta v_{i}+\omega_{2,i}, v_{i}) \coloneqq & F_2(s_{i}, \Delta v_{i}, v_{i}, \omega_{1,i}, \omega_{2,i}).   \label{eq2.10}
\end{eqnarray}

Considering the nature of false data injection attacks, they are assumed to be launched based on ACC sensor measurements perceived by the attacker~\citep{wang2023novel}; in other words, $\omega_{1,i}$ and $\omega_{2,i}$ are considered as function of the variable $s_i$ and $\Delta v_i$, respectively. Specifically, a class of such candidate attacks without being subject to any specific statistical distribution can be described by
\begin{eqnarray}
    \omega_{1,i} = h_1(s_i) = \delta_{1,i}s_{i}, ~ \omega_{2,i} = h_2(\Delta v_i) = \delta_{2,i}\Delta v_{i}.   \label{eq2.11}
\end{eqnarray}
Similar to $\delta_{i}$ of equation~\eqref{eq2.7}, $\delta_{1,i}$ and $\delta_{2,i}$ are also unstructured variables assumed to be bounded by some known bounds, e.g., $\delta_{j,i} \in \Upsilon_j \coloneqq \left\{\delta_{j,i}: \left|\delta_{j,i}\right| \leq z_j, z_j \in \mathbb{R}_{\geq 0}\right\}$, $j=1,2$, with $z_1=z_2=0$ indicating no attack. While it is not feasible to exhaustively characterize the analytical form of every attack possible and we do not intend to do so in this study, the above expressions of~\eqref{eq2.11} represent a broad class of candidate attacks based on which we are interested to derive their mathematical characterization. 

Considering~\eqref{eq2.11}, the mixed traffic dynamics of vehicle acceleration under Type~II attacks is given by
\begin{numcases}{\dot{v}_{i} =}
    f_{\text{HDV}}(s_{i}, \Delta v_{i}, v_{i}), ~\text{if}~ i \in \mathcal{H},  \label{eq2.12a}  \\
    f_{\text{ACC}}(s_{i}, \Delta v_{i}, v_{i}), ~\text{if}~ i \in \mathcal{A} ~\text{not attacked},    \label{eq2.12b}  \\ 
    F_2(s_{i}, \Delta v_{i}, v_{i}, \delta_{1,i}, \delta_{2,i}), ~\text{if}~ i \in \mathcal{A} ~\text{attacked}.   \label{eq2.12c}
\end{numcases}

\subsection{Illustration of mixed traffic under cyberattack}\label{section2.4}

Following the general framework presented above, we further illustrate the mixed traffic flow dynamics with concrete car-following models for ease of mathematical analysis and numerical study. For HDVs, the intelligent driver model
(IDM) is adopted since it has been shown to be able to accurately describe car-following dynamics of human drivers~\citep{treiber2000congested} and is also widely used in the traffic engineering community thanks to its favorable performance~\citep{talebpour2016influence}. The IDM is a multi-regime model which can provide a high level of realism in capturing the dynamics of different congestion levels~\citep{sarker2019review}. In addition, recent studies have shown that the IDM can well replicate human driving behavior with a high degree of accuracy and outperform other car-following models based on the examination of real-world driving data~\citep{pourabdollah2017calibration,he2023calibrating}.

For any HDV, $i \in \mathcal{H}$, following the IDM equation~\eqref{eq2.9a} is explicitly given by
\begin{eqnarray}
    &~&\hskip-20pt f_{\text{HDV}} = a \left[ 1 - \left(\frac{v_{i}}{v_{0}}\right)^{4}  - \left(\frac{s^{\ast}(v_{i},\Delta v_{i})}{s_{i}}\right)^{2} \right],   \label{eq2.13}
\end{eqnarray}
with
\begin{equation}
    s^{\ast}(v_{i},\Delta v_{i}) = s_{0} + \max\left\{0, v_{i}T - \frac{v_{i}\Delta v_{i}}{2\sqrt{ab}}\right\},    \label{eq2.14}
\end{equation}
where $a$ is the maximum acceleration, $b$ is the comfortable braking deceleration, $v_{0}$ is the desired speed, $s_{0}$ is the minimum spacing, $T$ is the desired time headway, and $l_{i-1}$ is the length of vehicle $i-1$. 

As shown in~\citep{ye2022car}, ACC vehicles tend to exhibit different driving behaviors from human drivers. While ACC vehicles made by distinct manufacturers may have different acceleration functions given by equation~\eqref{eq2.5b}, the optimal velocity with relative velocity (OVRV) model~\citep{milanes2013cooperative} is taken to describe ACC car-following dynamics for multiple reasons. The OVRV model has been widely used for ACC systems~\citep{milanes2013cooperative}. It follows a constant time-gap policy which is consistent with the practical implementation of intelligent vehicles~\citep{ioannou1993autonomous}. The OVRV model has also been shown to fit well to both simulated and real trajectories of vehicles equipped with ACC capabilities~\citep{milanes2013cooperative}.

Following the OVRV model, equation~\eqref{eq2.9b} is written as
\begin{eqnarray}
    f_{\text{ACC}} = k_{1}\left(s_{i} - \eta - \tau v_{i}\right) + k_{2}\Delta v_{i}, ~ i \in \mathcal{A}, ~ \text{not attacked},       \label{eq2.15}
\end{eqnarray}
where $\eta$ is the jam distance, i.e., inter-vehicle spacing at rest, $\tau$ is the desired time gap, $k_{1}$ and $k_{2}$ are positive gains on the effective time gap and the relative speed, respectively. While many car-following variables are time-variant, the argument $t$ is omitted for brevity wherever appropriate.

Consequently, following Section~\ref{section2.2} the acceleration dynamics $F_1$ resulting from Type~I attacks is given by
\begin{eqnarray}
    F_1 = k_{1}\left(s_{i} - \eta - \tau v_{i}\right) + k_{2}\Delta v_{i} + \delta_{i} v_{i}, ~ i \in \mathcal{A}, ~ \text{attacked}.   \label{eq2.16}
\end{eqnarray}
Similarly, it follows from Section~\ref{section2.3} that $F_2$ corresponding to Type~II attacks is 
\begin{eqnarray}
    F_2 = k_{1}\left(s_{i} + \delta_{1,i}s_{i} - \eta - \tau v_{i}\right) + k_{2}\left(\Delta v_{i} + \delta_{2,i}\Delta v_{i}\right), ~ i \in \mathcal{A}, ~ \text{attacked}.   \label{eq2.17}
\end{eqnarray}

\section{Analytical Characterization of Type~I Attacks on ACC Vehicles}\label{section3}

In this section, we first present some results on linear stability analysis of car-following dynamics. These results are then leveraged to synthesize a class of malicious attacks on ACC control commands, destabilizing ACC car-following dynamics. Further, we relax the conditions on attacks to cover a broader class of candidate attacks that can degrade, but may not necessarily destabilize, ACC dynamics. 

\subsection{Linear stability analysis}\label{section3.1}
Equation~\eqref{eq2.2}, or specifically equations~\eqref{eq2.5a} and~\eqref{eq2.5b}, describe vehicle driving behavior in terms of car-following dynamics. To ensure practical driving behavior of any vehicle being considered the following rational driving constraints~(RDC)~\citep{wilson2011car} are expected to be satisfied
\begin{eqnarray}
    \beta_{1} \coloneqq \frac{\partial \dot{v}}{\partial s} > 0, ~~ \beta_{2} \coloneqq \frac{\partial \dot{v}}{\partial \Delta v} > 0, ~~ \beta_{3} \coloneqq \frac{\partial \dot{v}}{\partial v} < 0. \label{eq3.1}
\end{eqnarray}
The expressions above can be understood as follows~\citep{wilson2011car}: $\beta_{1} > 0$ says that a larger spacing shall result in more acceleration; $\beta_{2} > 0$ indicates that a greater relative speed shall lead to more acceleration; $\beta_{3} < 0$ signifies that a vehicle tends to accelerate less as its speed gets larger. The RDC ensure a simple criterion for the existence of rational car-following models. Interested reader is referred to~\citep{wilson2011car} for a detailed interpretation on the expressions of~\eqref{eq3.1}.

The RDC ensure a simple criterion for the existence of car-following models. With the partial derivative of $\dot{v}$ with respect to $(s,\Delta v,v)$ evaluated at the equilibrium flow solution $(s_{\text{e}}, 0, v_{\text{e}})$, a car-following model is said to be string stable around the equilibrium state if the following $\lambda_{2}$ criteria hold~\citep{wilson2011car}
\begin{equation}
    \lambda_2 = \frac{\beta_{1}}{\beta_{3}^3}\left(\frac{\beta_{3}^2}{2}-\beta_{2}\beta_{3}-\beta_{1}\right) < 0.    \label{eq3.2}
\end{equation}

\begin{remark}\label{Remark3.1}
The widely used $\lambda_{2}$ criteria shown above are derived based on linear stability analysis considering an infinite traffic flow~\citep{wilson2011car}. They have been shown effective in evaluating string stability of both homogeneous and heterogeneous traffic flow described by car-following models~\citep{wilson2011car,talebpour2016influence}. The $\lambda_{2}$ condition follows the underlying assumption of an infinite traffic flow for theoretical soundness~\citep{montanino2021string}. In the present work, we follow this convention in the subsequent mathematical analysis and numerical study.
\end{remark}

\subsection{Synthesis of malicious Type~I attacks destabilizing traffic}\label{section3.2}

Following the results presented in Section~\ref{section2.2}, we analytically characterize the malicious nature of Type~I attacks based on the results given in Section~\ref{section3.1}. Specifically, we derive analytical properties of $\delta_{i}$ for attacks to destabilize ACC dynamics. To this end, we introduce the following physically interpretable conditions based on the understanding of vehicle driving behavior and attacker intentions. 
\begin{itemize}
    \item[(i)] Rational driving behavior described by the RDC of~\eqref{eq3.1} shall be maintained for the attacked ACC vehicle since irrational driving could easily get detected and identified. This essentially ensures that attacks are launched in a subtle manner. 
    \item[(ii)] Candidate attacks shall compromise the driving behavior of the attacked ACC vehicle due to their malicious nature. This will be analytically characterized using the $\lambda_2$ condition introduced in Section~\ref{section3.1}. 
\end{itemize}
As shown in~\citep{li2022detecting}, commonly used random attacks, like Gaussian noise, on ACC vehicles can be well detected using machine learning based anomaly detection techniques. It is also reasonable to assume that erratic, uncalculated and aggressive attacks can be effectively detected by intrusion detection systems~\citep{khraisat2019survey}. Hence, for a good degree of realism it is necessary to synthesize candidate attacks taking into account the driving behavior of vehicles and the possible intention of malicious actors, appearing to have been largely ignored in relevant studies. Essentially, the conditions (i) and (ii) indicate that attackers tend to alter the driving behavior of compromised ACC vehicles in a subtle manner while remaining malicious. In fact, even subtle changes to vehicle driving behavior could result in widespread disruption to the transportation network, causing substantial traffic congestion and excessive energy consumption and emissions from vehicles~\citep{dong2020impact,li2023exploring}. 

\begin{figure}[t!]
  \centering
  \subfloat[$-0.5 \leq \delta_{i} \leq 1.0$]{\includegraphics[width=0.45\textwidth]{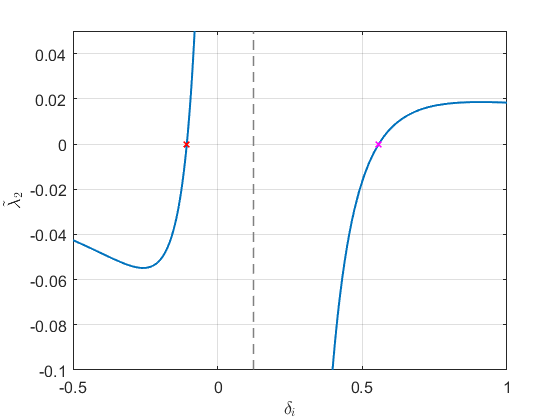}\label{lambda2-plot1}}
  \subfloat[$-4 \leq \delta_{i} \leq 5$]{\includegraphics[width=0.45\textwidth]{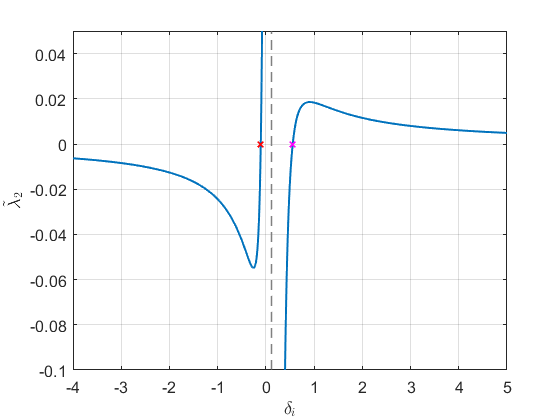}\label{lambda2-plot2}}
  \caption{Illustration of the function $\tilde{\lambda}_2$ shown in equation~\eqref{eq3.7} with different ranges of $\delta_{i}$, (a) $\delta_{i} \in \left[-0.5, 1.0\right]$ and (b) $\delta_{i} \in \left[-4, 5\right]$, where $k_1 = 0.05$, $k_2 = 0.10$, and $\tau = 2.50$. The vertical dashed line corresponds to $\delta_{i} = \tau k_1$. It is consistent with the analysis that $\tilde{\lambda}_2$ has two zeros highlighted with red markers, with the larger one being infeasible due to contradiction with~\eqref{eq3.5}. Clearly, for the portion of $\tilde{\lambda}_2$ on the left-hand side of the dashed line, i.e., $\delta_{i} < \tau k_1$, $\tilde{\lambda}_2 \geq 0$ when $\delta_{i} \geq \delta_{i}^*$.}\label{lambda2}%
\end{figure}

We first examine the mathematical implications of condition~(i) on candidate attacks shown in equation~\eqref{eq2.7}. With the attacked ACC dynamics given by equation~\eqref{eq2.9c}, to ensure the RDC it follows that
\begin{eqnarray}
    &~& \tilde{\beta}_{1} \coloneqq {\partial F_1}/{\partial s} = k_{1} > 0,  \label{eq3.3}  \\
    &~& \tilde{\beta}_{2} \coloneqq {\partial F_1}/{\partial \Delta v} = k_{2} > 0,    \label{eq3.4}  \\ 
    &~& \tilde{\beta}_{3} \coloneqq {\partial F_1}/{\partial v} = -\tau k_{1} + \delta_{i} < 0.  \label{eq3.5}
\end{eqnarray}
Since the expression~\eqref{eq3.5} needs to hold for any attack $\delta_{i} \in \Lambda$, it follows that 
\begin{eqnarray}
    -\tau k_{1} + r < 0 \Longrightarrow r < \tau k_{1}.   \label{eq3.6}
\end{eqnarray}
Plugging the expressions~\eqref{eq3.3},~\eqref{eq3.4} and~\eqref{eq3.5} into~\eqref{eq3.2} yields
\begin{eqnarray}
    \tilde{\lambda}_2 = {\scriptstyle \frac{k_{1}}{\left(-\tau k_{1} + \delta_{i}\right)^3}\left[\frac{\left(-\tau k_{1} + \delta_{i}\right)^2}{2}-k_{2}\left(-\tau k_{1} + \delta_{i}\right)-k_{1}\right]} = {\scriptstyle \frac{-k_{1}^2}{\left(-\tau k_{1} + \delta_{i}\right)^3} - \frac{k_{1}k_{2}}{\left(-\tau k_{1} + \delta_{i}\right)^2} + \frac{k_{1}}{2\left(-\tau k_{1} + \delta_{i}\right)}}.    \label{eq3.7}
\end{eqnarray}
For a candidate attack to destabilize the ACC dynamics, it is necessary to have $\tilde{\lambda}_2 > 0$. We shall first derive the solutions to the equation $\tilde{\lambda}_2 = 0$. Following equation~\eqref{eq3.7} with basic algebraic calculation, one can verify that the solutions of $\tilde{\lambda}_2 = 0$ are given by
\begin{eqnarray}
    \delta_{i} = \tau k_{1} + k_{2} + \sqrt{k_{2}^2+2k_{1}} ~\text{or}~ \tau k_{1} + k_{2} - \sqrt{k_{2}^2+2k_{1}}.       \label{eq3.8}
\end{eqnarray}
It it easy to determine that $\delta_{i} = \tau k_{1} + k_{2} + \sqrt{k_{2}^2+2k_{1}} > \tau k_{1}$ is not feasible since it contradicts the fact that $\delta_{i} \leq r < \tau k_{1}$ according to~\eqref{eq3.6}. For the other solution $\delta_{i} = \tau k_{1} + k_{2} - \sqrt{k_{2}^2+2k_{1}}$ to be feasible, it follows from $\left|\delta_{i}\right| \leq r$ that
\begin{eqnarray}
    -r \leq \delta_{i} = \tau k_{1} + k_{2} - \sqrt{k_{2}^2+2k_{1}} \leq r,   \label{eq3.9}
\end{eqnarray}
leading to 
\begin{subequations}
\begin{equation}
    r \geq -r_{1} \coloneqq -\tau k_{1} - k_{2} + \sqrt{k_{2}^2+2k_{1}},  \label{eq3.10a}
\end{equation}
\begin{equation}
    r \geq r_{1} \coloneqq \tau k_{1} + k_{2} - \sqrt{k_{2}^2+2k_{1}}.   \label{eq3.10b}
\end{equation}
\end{subequations}
Applying the condition of $r < \tau k_{1}$ shown in~\eqref{eq3.6} to~\eqref{eq3.10a} yields 
\begin{eqnarray}
    \tau k_{1} > -\tau k_{1} - k_{2} + \sqrt{k_{2}^2+2k_{1}},   \label{eq3.11}
\end{eqnarray}
resulting in
\begin{eqnarray}
    2\tau^2 k_{1} + 2\tau k_{2} - 1 > 0.   \label{eq3.12}
\end{eqnarray}
Further, it is easy to verify that the conditions~\eqref{eq3.6} and~\eqref{eq3.10b} can always be satisfied simultaneously. Hence, for the existence of a feasible solution to the equation $\tilde{\lambda}_{2} = 0$ with the RDC satisfied, it is necessary that 
\begin{eqnarray}
   \tau k_{1} > r \geq \max\left\{-r_{1}, r_{1}, 0\right\},   \label{eq3.13}
\end{eqnarray}
with $r_{1}$ defined in~\eqref{eq3.10b}, and that the following condition on model parameters holds
\begin{eqnarray}
   2\tau^2 k_{1} + 2\tau k_{2} - 1 > 0.  \label{eq3.14}
\end{eqnarray}
For $\delta_{i} \in \Lambda$, the function $\tilde{\lambda}_{2}$ has one zero with~\eqref{eq3.13} and~\eqref{eq3.14} satisfied, i.e., $\delta_{i}^* = \tau k_{1} + k_{2} - \sqrt{k_{2}^2+2k_{1}}$, and it continues to increase as $\delta_{i}$ becomes larger than $\delta_{i}^*$, with a numerical example shown in Fig.~\ref{lambda2}.

By virtue of~\eqref{eq3.2}, for a candidate attack to destabilize an ACC model, it is necessary that $\tilde{\lambda}_2 \geq 0$. Hence, it follows that
\begin{eqnarray}
   \delta_{i}^* \leq \delta_{i} \leq r.  \label{eq3.15}
\end{eqnarray}
Letting $\Omega$ denote a class of Type~I attacks destabilizing the ACC dynamics, it follows from~\eqref{eq3.13},~\eqref{eq3.14} and~\eqref{eq3.15} that
\begin{eqnarray}
    \Omega \coloneqq \bigl\{ \omega_{i}: \delta_{i} \in \left[\delta_{i}^*, r\right], ~\text{with}~ \tau k_{1} > r \geq \max\left\{-r_{1}, r_{1}, 0\right\} ~\text{and}~ 2\tau^2 k_{1} + 2\tau k_{2} - 1 > 0 \bigr\},    \label{eq3.16}
\end{eqnarray}
where $r_{1}$ is defined in~\eqref{eq3.10b}, and $k_{1}$ and $k_{2}$ are positive gains of the OVRV model.

\begin{remark}\label{remark3.2}
While we do not intend to exhaustively characterize the analytical form of all possible attacks, the set $\Omega$ represents a broad class of malicious candidate attacks on ACC control commands. Although attacks could be launched also from the complement of $\Omega$, they may result in noticeably irrational driving behavior of the compromised vehicle, likely making it easier to be detected and identified. 
\end{remark}

\subsection{Relaxed conditions on Type~I attacks degrading traffic}\label{section3.3}

In the previous subsection, we have derived analytical conditions, given by equation~\eqref{eq3.16}, on a set of candidate attacks that can destabilize ACC dynamics in the context of linear stability analysis. While destabilizing traffic could be well aligned with the malicious nature of attackers, it may be demanding to achieve such a goal as the space $\Omega$ could be limited. To this end, we relax the conditions derived in Section~\ref{section3.2} to synthesize a class of candidate attacks that can degrade, but may not necessarily destabilize, ACC dynamics. 

In the absence of attack, it follows from~\eqref{eq3.2} that the value of $\lambda_{2}$ corresponding to an ACC vehicle is given by
\begin{equation}
    \lambda_2 = \frac{k_{1}}{(-\tau k_{1})^3}\left(\frac{(-\tau k_{1})^2}{2}-k_{2}(-\tau k_{1})-k_{1}\right).    \label{eq3.17}
\end{equation}
The counterpart of $\lambda_{2}$, in the presence of an attack, is given by $\tilde{\lambda}_{2}$ in equation~\eqref{eq3.7}. For a candidate attack to degrade ACC stability, it is necessary to have
\begin{equation}
    \tilde{\lambda}_{2} > \lambda_2.    \label{eq3.18}
\end{equation}
In other words, one has
\begin{eqnarray}
    \frac{-k_{1}^2}{\left(-\tau k_{1} + \delta_{i}\right)^3} - \frac{k_{1}k_{2}}{\left(-\tau k_{1} + \delta_{i}\right)^2} + \frac{k_{1}}{2\left(-\tau k_{1} + \delta_{i}\right)} > \frac{-k_{1}^2}{\left(-\tau k_{1}\right)^3} - \frac{k_{1}k_{2}}{\left(-\tau k_{1}\right)^2} + \frac{k_{1}}{2\left(-\tau k_{1}\right)},    \label{eq3.19}
\end{eqnarray}
which is equivalent to 
\begin{eqnarray}
    p(\delta_i) \coloneqq (\tau^2k_1 + 2\tau k_2 - 2)(\delta_i - \tau k_1)^3 + \tau^3k_1^2(\delta_i - \tau k_1)^2 - 2\tau^3k_1^2k_2(\delta_i - \tau k_1) - 2\tau^3 k_1^3 < 0.    \label{eq3.20}
\end{eqnarray}
Hence, for a candidate attack to degrade, but not necessarily destabilize, ACC dynamics, i.e., $\tilde{\lambda}_{2} > \lambda_2$, the following conditions shall be satisfied
\begin{eqnarray}
    &~& p(\delta_i) < 0,  \label{eq3.21}  \\
    &~& \delta_i \in \Lambda,    \label{eq3.22}  \\ 
    &~& r < \tau k_{1}.  \label{eq3.23}
\end{eqnarray}
where~\eqref{eq3.23} is same as~\eqref{eq3.6} for ensuring the RDC. 

\begin{remark}\label{remark3.3}
The conditions on attacks derived in Section~\ref{section3.3} are more relaxed than those of Section~\ref{section3.2}, considering $\lambda_2 < 0$. However, ACC dynamics may be unstable in the first place (without attack), i.e., $\lambda_2 > 0$. Consequently, the conditions given by~\eqref{eq3.21}--\eqref{eq3.23} become stronger than those of equation~\eqref{eq3.16}, i.e., $\tilde{\lambda}_{2} > \lambda_{2} > 0$. For implementation, one may simply apply the conditions obtained in Section~\ref{section3.2} for attack synthesis. This will be further illustrated via numerical study in Section~\ref{section5}.
\end{remark}

\section{Analytical Characterization of Type~II Attacks on ACC Vehicles}\label{section4}

In this section, we leverage the stability results presented in Section~\ref{section3.1} to synthesize a class of malicious attacks on ACC sensor measurement, destabilizing ACC dynamics. Further, we relax the conditions on attacks to cover a class of candidate attacks that can degrade, but may not necessarily destabilize, ACC dynamics. 

\subsection{Synthesis of malicious Type~II attacks destabilizing traffic}\label{section4.1}

\begin{figure}[t!]
  \centering
  \subfloat[Value of $\hat{\lambda}_2$ by~\eqref{eq4.5}]{\includegraphics[width=0.45\textwidth]{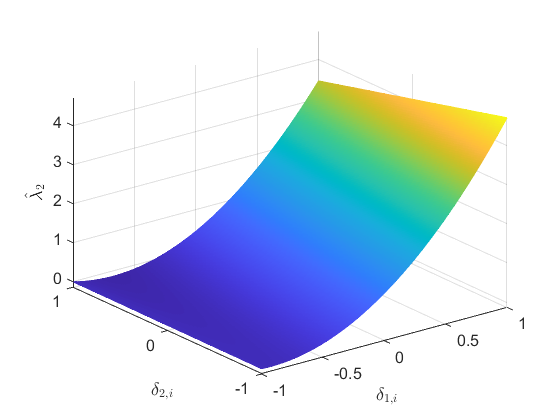}\label{lambda2-hat}}
  \subfloat[Value of $\theta$ by~\eqref{eq4.16}]{\includegraphics[width=0.45\textwidth]{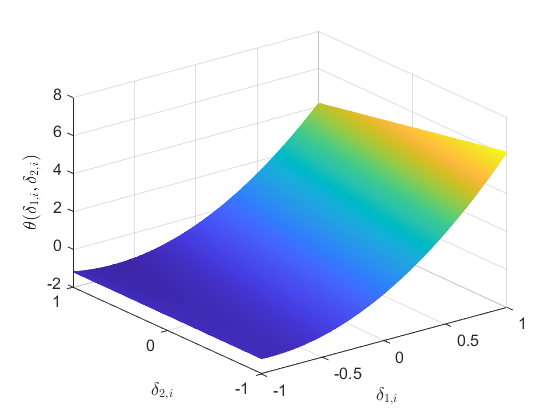}\label{theta_delta}}
  \caption{A graphic illustration of $\hat{\lambda}_2$ given by~\eqref{eq4.5} and $\theta$ shown in~\eqref{eq4.16} with $-1 \leq \delta_{1,i}, \delta_{2,i} \leq 1$, where $k_1 = 0.05$, $k_2 = 0.10$, and $\tau = 2.50$.}\label{lambda2_hat}%
\end{figure}

Similar to Section~\ref{section3.2}, we first examine the mathematical implications of condition~(i) on candidate attacks shown in~\eqref{eq2.11}. With the attacked ACC dynamics given by equation~\eqref{eq2.12c}, to ensure the RDC it follows that
\begin{eqnarray}
    &~& \hat{\beta}_{1} \coloneqq {\partial F_2}/{\partial s} = k_{1}(1+\delta_{1,i}) > 0,  \label{eq4.1}  \\
    &~& \hat{\beta}_{2} \coloneqq {\partial F_2}/{\partial \Delta v} = k_{2}(1+\delta_{2,i}) > 0,    \label{eq4.2}  \\ 
    &~& \hat{\beta}_{3} \coloneqq {\partial F_2}/{\partial v} = -\tau k_{1} < 0.  \label{eq4.3}
\end{eqnarray}
It follows from~\eqref{eq4.1} and~\eqref{eq4.2} that 
\begin{eqnarray}
    1 + \delta_{1,i} > 0 ~\text{and}~ 1 + \delta_{2,i} > 0.   \label{eq4.4}
\end{eqnarray}
Plugging the expressions~\eqref{eq4.1},~\eqref{eq4.2} and~\eqref{eq4.3} into~\eqref{eq3.2} yields
\begin{eqnarray}
    \hat{\lambda}_2 = {\scriptstyle -\frac{k_{1}\left(1+\delta_{1,i}\right)}{\left(\tau k_{1}\right)^3}\bigl[\frac{\left(\tau k_{1}\right)^2}{2}+\tau k_{1}k_{2}\left(1+\delta_{2,i}\right)-k_{1}\left(1+\delta_{1,i}\right)\bigr]}.   \label{eq4.5}
\end{eqnarray}
Clearly, $\hat{\lambda}_2$ is a function of both $\delta_{1,i}$ and $\delta_{2,i}$, with a graphic illustration shown in Fig.~\ref{lambda2-hat}. In order for a candidate attack to destabilize the ACC model, it is necessary to have $\hat{\lambda}_2 > 0$. While the sign of $\partial \hat{\lambda}_2/\partial \delta_{1,i}$ is indefinite, it is interesting to note that
\begin{eqnarray}
    \frac{\partial \hat{\lambda}_2}{\partial \delta_{2,i}} = -\frac{k_{2}(1+\delta_{1,i})}{\tau^2k_{1}} < 0,     \label{eq4.6}
\end{eqnarray}
indicating that the value of $\hat{\lambda}_2$ monotonically decreases with the increase of $\delta_{2,i}$. Hence, one may tend to launch attacks with a smaller value of $\delta_{2,i}$ to destabilize ACC dynamics. Letting $\Theta$ denote a set of Type~II attacks destabilizing the ACC model, it directly follows that
\begin{eqnarray}
    \Theta \coloneqq \left\{ \omega_{j,i}: \hat{\lambda}_2 > 0, ~\delta_{j,i} > -1, ~\delta_{j,i} \in \Upsilon_j, ~j=1,2 \right\},    \label{eq4.14}
\end{eqnarray}
with $\hat{\lambda}_2$ given by equation~\eqref{eq4.5}.

\subsection{Relaxed conditions on Type~II attacks degrading traffic}\label{section4.2}

As mentioned before, while destabilizing traffic is well aligned with the malicious nature of potential attacks, e.g., given by equation~\eqref{eq4.14}, it may be demanding to achieve this goal as the attack space $\Theta$ could be limited. Hence, we relax the conditions derived in Section~\ref{section4.1} to synthesize a class of candidate attacks that can degrade, but may not necessarily destabilize, ACC dynamics. 

Similar to~\eqref{eq3.18}, for a Type~II attack to degrade ACC stability, it is necessary to have
\begin{equation}
    \hat{\lambda}_{2} > \lambda_2,    \label{eq4.15}
\end{equation}
which, after some algebraic manipulation, leads to
\begin{eqnarray}
    \theta(\delta_{1,i},\delta_{2,i}) \coloneqq 2\delta_{1,i}^2 - (\tau^2 k_1+2\tau k_2-4)\delta_{1,i} - 2\tau k_2\delta_{2,i} - 2\tau k_2\delta_{1,i}\delta_{2,i} > 0.    \label{eq4.16}
\end{eqnarray}
It is easy to verify that
\begin{eqnarray}
    \frac{\partial \theta}{\partial \delta_{2,i}} = -2\tau k_2(1+\delta_{1,i}) < 0,     \label{eq4.17}
\end{eqnarray}
indicating that the value of $\theta$ monotonically decreases with the increase of $\delta_{2,i}$. Hence, one may tend to launch attacks with a smaller value of $\delta_{2,i}$ to degrade ACC dynamics, consistent with the indication of~\eqref{eq4.6}. A graphic illustration of the function $\theta$ is presented in Fig.~\ref{theta_delta}. However, it is noted that the sign of $\partial \theta/\partial \delta_{2,i}$ is indefinite. Thus, for a Type~II attack to degrade ACC stability without necessarily having to destabilize ACC dynamics, i.e., $\hat{\lambda}_{2} > \lambda_2$, the following conditions shall be satisfied
\begin{eqnarray}
    &~& \theta(\delta_{1,i},\delta_{2,i}) > 0,  \label{eq4.18}  \\
    &~& \delta_{j,i} \in \Upsilon_j, ~j=1,2,    \label{eq4.19}  \\ 
    &~& 1 + \delta_{j,i} > 0, ~j=1,2.  \label{eq4.20}
\end{eqnarray}
where~\eqref{eq4.20} is same as~\eqref{eq4.4} for ensuring the RDC. 

\begin{remark}\label{remark4.1}
While Type~I and Type~II attacks could occur to an ACC vehicle simultaneously, it is straightforward to carry out a similar analysis as seen in Sections~\ref{section3} and~\ref{section4} by considering the compromised ACC dynamics involving both equations~\eqref{eq2.7} and~\eqref{eq2.11}. This is omitted for brevity. 
\end{remark}

\section{Numerical Results}\label{section5}

In this section, we present extensive numerical results to illustrate the mechanism of Type~I and~II attacks, along with their impacts on microscopic car-following dynamics and macroscopic traffic flow.

\subsection{Simulation settings}\label{section5.1}

\begin{figure}[t!]
	\centering
	\includegraphics[trim={7cm 5cm 1cm 5.2cm},clip, width=0.6\textwidth]{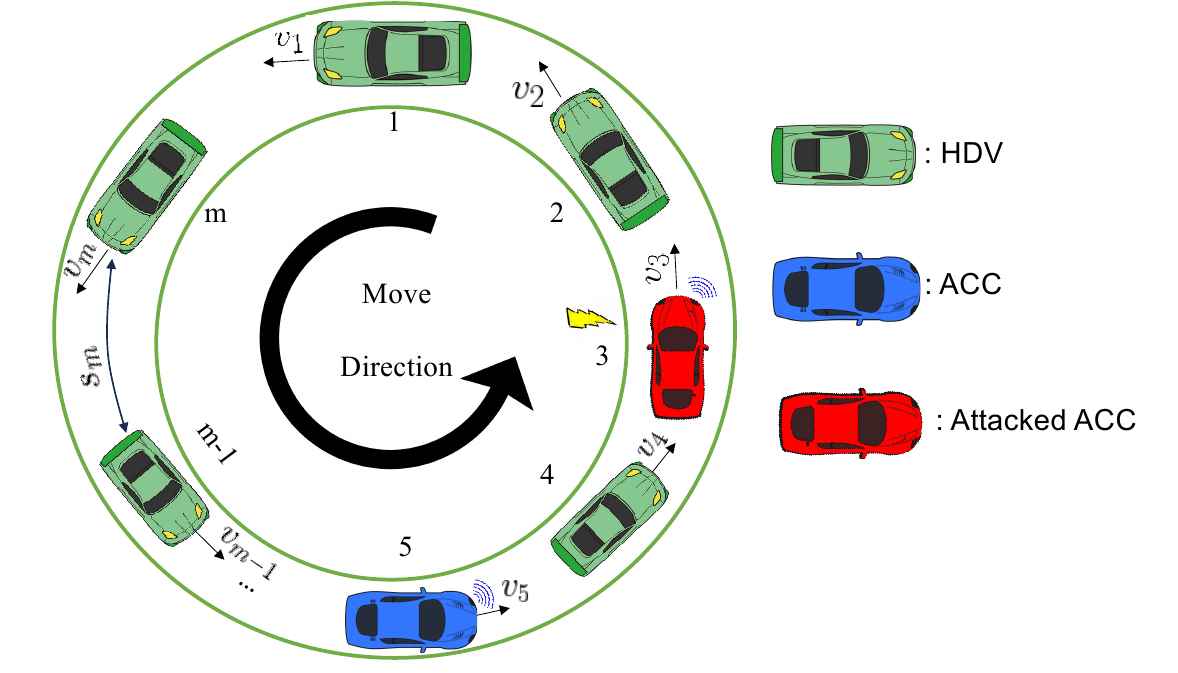}
	\caption{A graphic illustration of a ring of $m$ vehicles in mixed traffic with ACC vehicles under potential attack.}
	\label{Simulation_setting}
\end{figure}

\begin{table}[t!]
\setlength{\tabcolsep}{1pt}  
\caption{Parameter Values of the IDM and OVRV Model}\label{Table_parameters}
\begin{center}
 \begin{tabular}{c c c c}
 \hline \hline
 \textbf{Parameter} ~&~ \textbf{Description} ~&~ \textbf{IDM} ~&~ \textbf{OVRV}  \\  [0.5ex]
 \hline
 $v_{0}$ ~&~ desired speed (\text{m/s}) ~&~ 30.0 ~&~ - \\ [0.3ex]
 \hline
 $T$ ~&~ time gap (\text{s}) ~&~ 1.5 ~&~ - \\  [0.3ex]
 \hline
 $s_{0}$ ~&~ minimum spacing (\text{m}) ~&~ 2.0 ~&~ -\\[0.3ex]
 \hline
 $a$ ~&~ maximum acceleration (\text{m/s\textsuperscript{2}}) ~&~ 1.4 ~&~ -  \\  [0.3ex]
 \hline
 $b$ ~&~ comfortable deceleration (\text{m/s\textsuperscript{2}}) ~&~ 2.0 ~&~ - \\  [0.3ex]
 \hline
 $l_{i}$ ~&~ vehicle length (\text{m}) ~&~ 5.0 ~&~ 5.0  \\ [0.3ex]
 \hline
 $k_{1}$ ~&~ ACC gain parameter (\text{s\textsuperscript{-2}}) ~&~ - ~&~ 0.05   \\  [0.3ex]
 \hline
 $k_{2}$ ~&~ ACC gain parameter (\text{s\textsuperscript{-1}}) ~&~ - ~&~ 0.10 \\  [0.3ex]
 \hline
 $\eta$ ~&~ jam distance (\text{m}) ~&~ - ~&~ 21.51  \\  [0.3ex]
 \hline
 $\tau$ ~&~ desired time gap (\text{s}) ~&~ - ~&~ 2.50 \\ [0.3ex]
 \hline \hline
\end{tabular}
\end{center}
\end{table}

As shown in Fig.~\ref{Simulation_setting}, a ring track is utilized in simulation, for multiple reasons: (1) it is capable of creating traffic jams in the absence of a physical bottleneck while considering only longitudinal vehicle dynamics; (2) it allows for study of traffic flow while neglecting the amount of incoming traffic; (3) it has been widely studied in simulation and field experiments to capture phantom traffic jams~\citep{sugiyama2008traffic,wu2019tracking}.

In this study, we simulate a total of $m$ = 40 vehicles on a ring track with length $L=$~1,400~m. At the initial point, all vehicles start with the same spacing gap, i.e., bumper-to-bumper distance, of 30~m. For illustrative purposes, 20\% of the simulated traffic is considered ACC vehicles evenly distributed in the flow, while other penetration rates could also be studied. It is assumed that half of the ACC vehicles are attacked. The speed limit of simulation is set as 30~m/s, i.e., 108~km/h, to reflect real-world conditions. As introduced in Section~\ref{section2.4}, the IDM and OVRV model are used to describe HDV and ACC dynamics, respectively. Their parameter values are shown in Table~\ref{Table_parameters}, with the ones of the IDM taken from~\citep{treiber2013traffic} and those of the OVRV model adopted from~\citep{gunter2019modeling} with modifications for illustrative purposes.

\subsection{Type~I attacks}\label{section5.2}

For Type~I attacks, the value of $\delta_i$ is taken from $[-0.10, 0.12]$, ensuring the conditions shown in~\eqref{eq3.16}. We simulate traffic flow with the candidate attacks given by~\eqref{eq2.7} to assess their impacts on traffic oscillation, safety, and efficiency.

\subsubsection{Trajectories and safety impacts}

\begin{figure}[t!]
  \centering
  \subfloat[$\delta_i = -0.10$]{\includegraphics[width=0.5\textwidth]{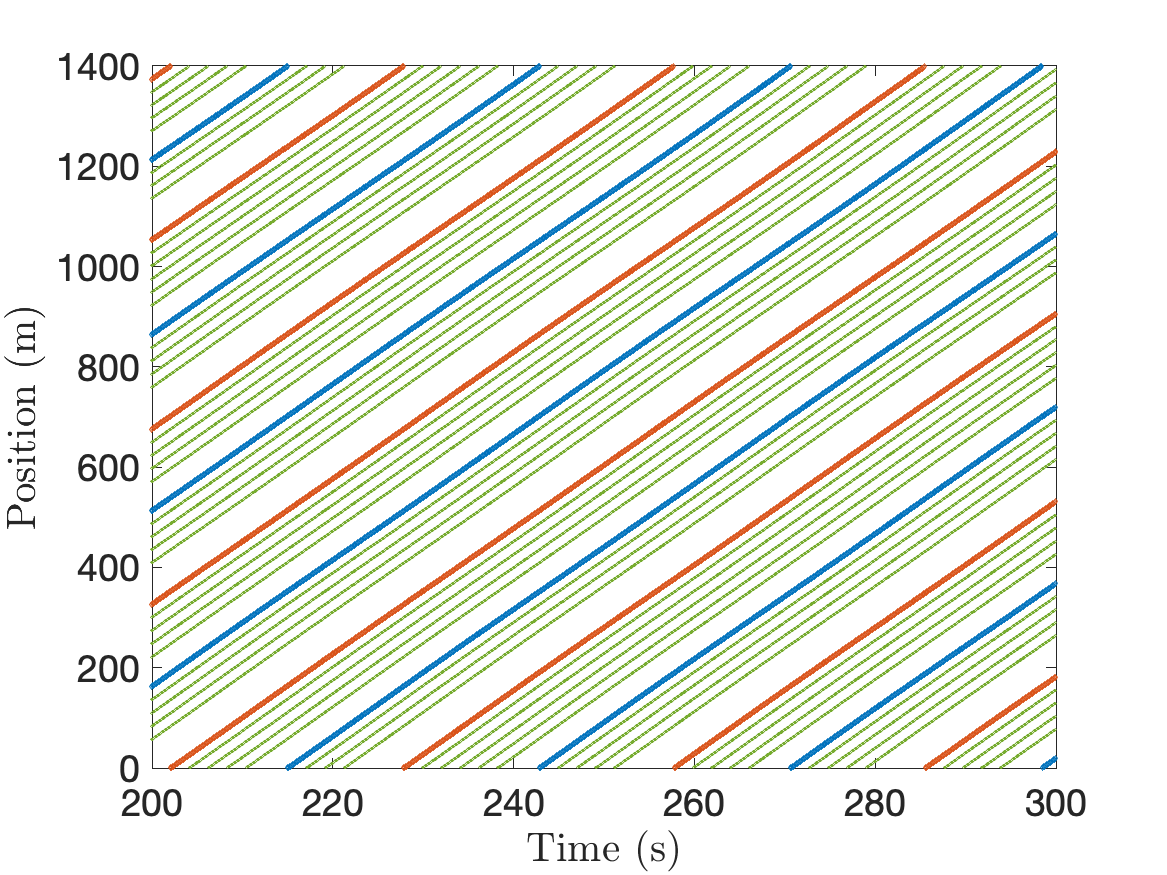}\label{trajectory_type_I_01}}
  \subfloat[$\delta_i = 0.00$]{\includegraphics[width=0.5\textwidth]{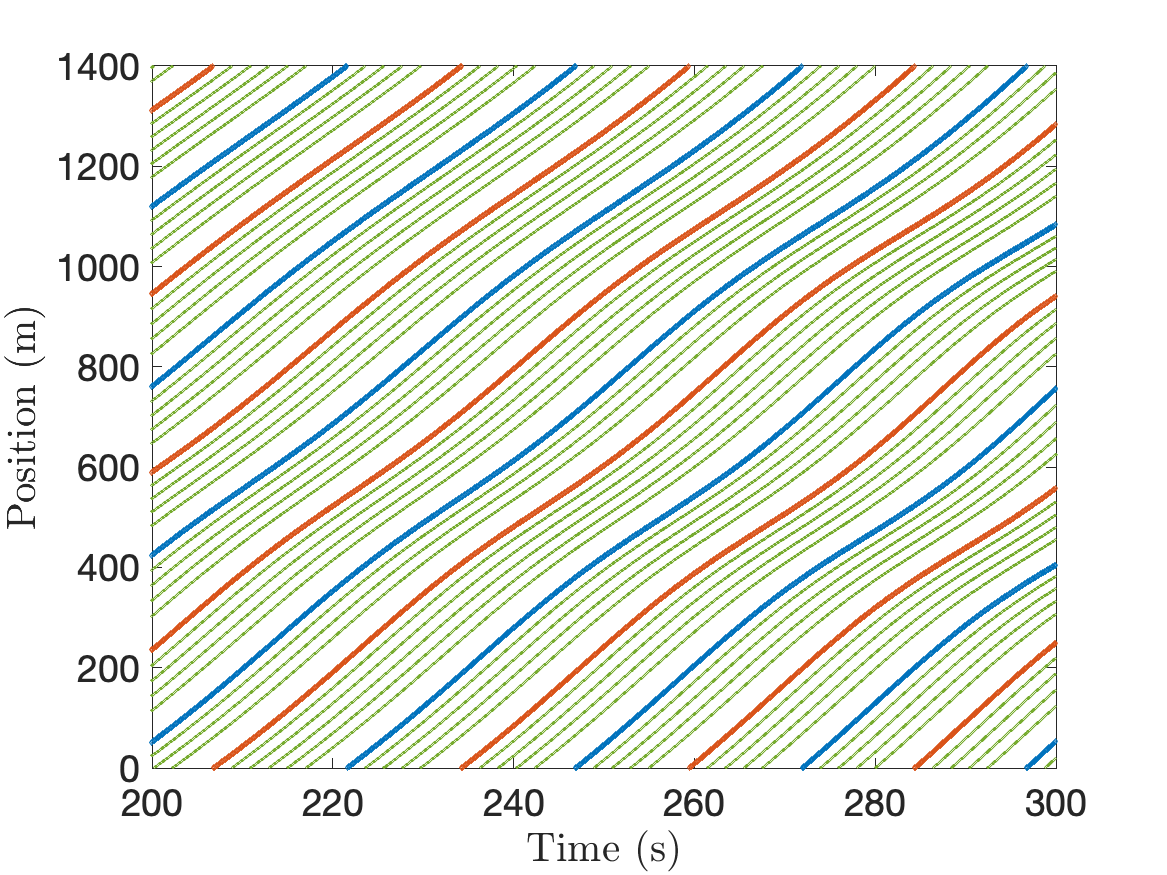}\label{trajectory_type_I_02}}
  \vfil%
  \vspace*{-1.0em}%
  \subfloat[$\delta_i = 0.06$]{\includegraphics[width=0.5\textwidth]{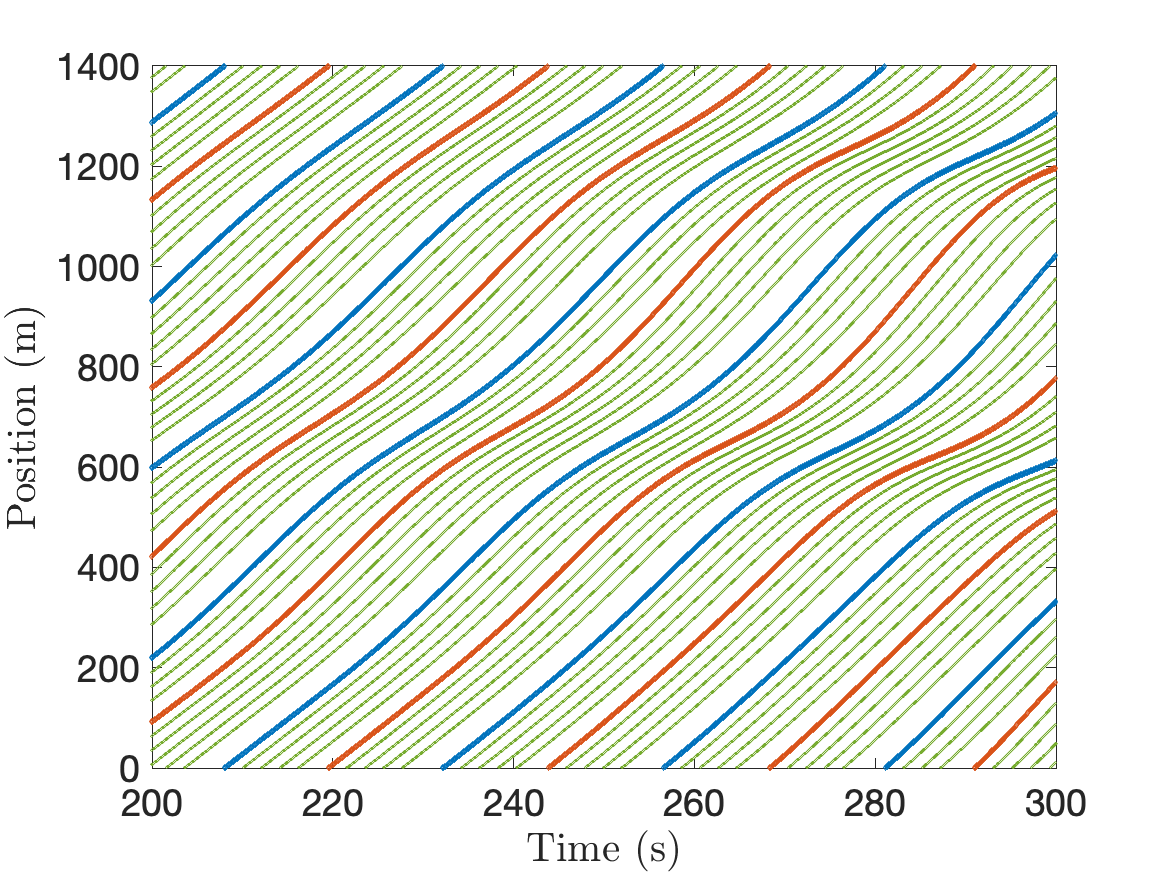}\label{trajectory_type_I_03}}
  \subfloat[$\delta_i = 0.10$]{\includegraphics[width=0.5\textwidth]{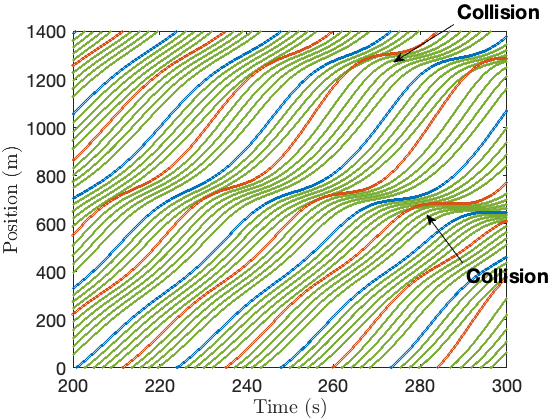}\label{trajectory_type_I_04}}
  \caption{Vehicle trajectory for different values of $\delta_i$. Greater traffic oscillations are observed when increasing the value of $\delta_i$, eventually resulting in rear-end collisions. Green, blue, and red correspond to HDV, normal ACC vehicle, and attacked ACC vehicle, respectively.}\label{trajectory_type_I}%
\end{figure}

Vehicle trajectories are presented in Fig.~\ref{trajectory_type_I} for four different values of $\delta_i$, with $\delta_i=0$ indicating no attack. It is observed that more traffic oscillations arise as the value of $\delta_i$ increases. Specifically, in Fig.~\ref{trajectory_type_I_01} with $\delta_i=-0.1$, less traffic oscillations occur to the attacked ACC vehicles compared with the non-attack scenario shown in Fig.~\ref{trajectory_type_I_02}. This is because a negative $\delta_i$, thereby a negative $\omega_{i}$ given by~\eqref{eq2.7}, incurred to the acceleration, would reduce speed of the attacked vehicle, resulting in a larger spacing gap between the compromised ACC vehicle and its lead vehicle.

\begin{table}[t!]
\setlength{\tabcolsep}{2.5pt}  
\caption{Ratio of dangerous risks ($p_r$) at different Time-to-Collision (TTC) thresholds}\label{results_TTC}
\begin{center}
 \begin{tabular}{c| c c c c c c }
 \hline
 \multirow{2}{*}{$\mathbf{\delta_{i}}$} & \multicolumn{6}{c}{\textbf{TTC}}\\ 
 \cline{2-7}
 ~&~ \textbf{1.5s}  ~&~ \textbf{2.0s} ~&~ \textbf{2.5s} ~&~ \textbf{3.0s} ~&~ \textbf{3.5s} ~&~ \textbf{4.0s} \\  [0.5ex]
 \hline
 0.00 ~&~ 0  ~&~ 0 ~&~ 0 ~&~ 0  ~&~ 0 ~&~ 0  \\  [0.3ex]
 \hline
 0.02 ~&~ 0 ~&~ 0 ~&~ 0 ~&~ 0  ~&~ 0 ~&~ 0 \\  [0.3ex]
 \hline
 0.04 ~&~ 0 ~&~ 0 ~&~ 0 ~&~ 0  ~&~ 0 ~&~ 0  \\  [0.3ex]
 \hline
 0.06 ~&~ 0 ~&~ 0 ~&~ 0 ~&~ 0  ~&~ 0 ~&~ 0.01\% \\ [0.3ex]
 \hline
 0.08 ~&~ 0.03\% ~&~ 0.04\%  ~&~ 0.04\%  ~&~ 0.05\%  ~&~ 0.05\%  ~&~  0.06\%   \\ [0.3ex]
 \hline 
 0.10 ~&~ 0.05\% ~&~ 0.60\% ~&~ 0.80\% ~&~ 1.00\% ~&~ 2.00\% ~&~ 2.00\% \\ [0.3ex]
 \hline 
 0.12 ~&~ 0.09\% ~&~ 1.00\% ~&~ 1.00\% ~&~ 2.00\%  ~&~ 3.00\% ~&~ 4.00\%  \\ [0.3ex]
 \hline 
\end{tabular}
\end{center}
\end{table}

When increasing $\delta_i$ to a positive value shown in Figs.~\ref{trajectory_type_I_03} and~\ref{trajectory_type_I_04}, greater traffic waves are observed due to the attacks, e.g., comparing Fig.~\ref{trajectory_type_I_03} to Fig.~\ref{trajectory_type_I_02}. However, increasing $\delta_i$ to 0.10 leads to rear-end collisions. Clearly, a positive $\delta_i$ tends to increase the acceleration of the attacked ACC vehicle, resulting in higher speed and insufficient spacing gap between vehicles. Consequently, a large, positive attack may lead to aggressive driving behavior of the compromised ACC vehicle.

Since aggressive driving could have significant impact on traffic safety, we compute the time-to-collision (TTC)~\citep{hayward1972near}, widely used in assessing safety-related car-following behavior~\citep{li2020analysis}, for safety analysis. Specifically, the TTC between vehicle $i$ and vehicle $i-1$ is calculated as
\begin{eqnarray}
    \text{TTC}_{i} = \begin{cases}
       \frac{s_i}{-\Delta v_i},  ~~~\text{if}~ \Delta v_i < 0 \label{eq: Delta_v <0},\\
       \infty,   ~~~~~~~\text{if}~ \Delta v_i \geq 0 \label{eq: Delta_v>0}.
           \end{cases}
\end{eqnarray}
A common approach for evaluating TTC is by computing the ratio of cases in dangerous risks (i.e., below the TTC safety threshold)~\citep{mahmud2019micro}, which is given by
\begin{equation}
    p_r = \frac{c_r}{c_\text{total}}\times 100\%,
\end{equation}
where $c_r$ is the number of cases in dangerous risks and $c_\text{total}$ is the total number of simulated cases. The ratio $p_r$ provides valuable insights into the safety performance of the traffic flow under consideration. A smaller $p_r$ indicates better traffic safety, while a larger one signifies a higher likelihood of dangerous situations where the TTC falls below the safety threshold. We set the threshold from 1.5~s to 4.0~s with increments of 0.5~s~\citep{mahmud2019micro}. 

Since a positive $\delta_i$ may result in rear-end collision, we consider its values ranging from 0.00 to 0.12 with increments of 0.02. The corresponding results are summarized in Table~\ref{results_TTC}. For $\delta_i =$ 0.00, 0.02, and 0.04, the $p_r$ is 0 for each TTC threshold, indicating a minimal chance of collision. However, when increasing the value of $\delta_i$ to 0.12, the value of $p_r$ is observed to increase to $0.09\%$ and $4\%$ for a TTC threshold of 1.5~s and 4~s, respectively. These results indicate a higher likelihood of collision when the attack value $\delta_i$ increases, consistent with the observations in Fig.~\ref{trajectory_type_I}.

\begin{figure}[t!]
  \centering
  \includegraphics[width=0.6\textwidth]{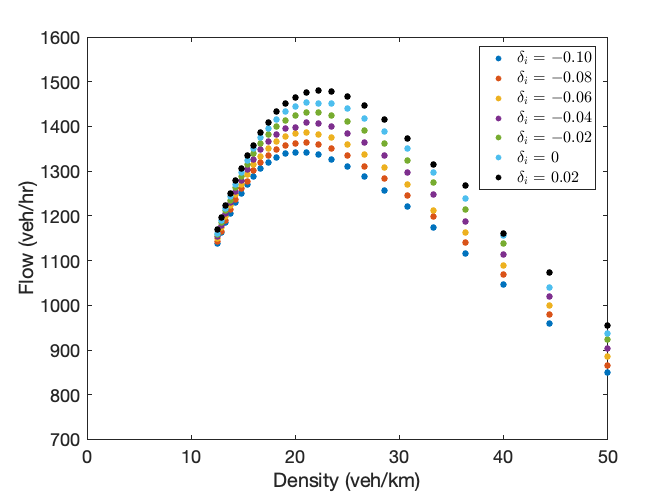}
  \caption{Fundamental diagrams of traffic flow at different values of $\delta_i$}
  \label{FD_type_I}
\end{figure}

\subsubsection{Traffic efficiency (fundamental diagram)}\label{sec: FDresultsI}

We have evaluated the microscopic impacts of Type~I attacks in the previous subsection. Here we are interested to study efficiency of the traffic flow from a macroscopic perspective. Studying the relationship between traffic flow and traffic density, commonly referred to as the fundamental diagram (FD), has been of interest to the transportation community for years~\citep{greenshields1934photographic}. Numerous studies have contributed to the understanding of the FD by observing and collecting traffic flow data~\citep{greenshields1934photographic,coifman2015empirical}. 

The FD characterizes the relationship between traffic density and traffic flow (flux). To establish this diagram, we conduct multiple simulations to collect the characteristic data of traffic flow. Using the simulation setup shown in Fig.~\ref{Simulation_setting}, we alter the loop length $L$, which in turn, changes the traffic density (denoted by $\rho$). During the simulation, we collect vehicle speeds to calculate the average speed ($\bar{v}$). The traffic flow rate (denoted by $q$) is obtained through the relationship, $q = \rho \bar{v}$. Hence, we obtain the FD, which displays the flow-density relationship, by varying the loop length. 

Fig.~\ref{FD_type_I} shows the FD with $\delta_i$ ranging from -0.10 to 0.02 with increments of 0.02. For a negative $\delta_i$, the resulting capacity is lower than that of the non-attack scenario. This is because a negative attack reduces vehicle speed and causes a larger spacing gap, resulting in degradation of traffic efficiency. In contrast, the capacity increases to some extent when $\delta_i$ is positive. However, it is worth noting that safety plays a more significant role for traffic flow in the presence of positive attacks. The potential risk of collision, along with the increased traffic oscillations, outweighs the marginal gain of capacity, making positive attacks less desirable in practical scenarios when designing candidate attacks that do not cause collision is prioritized~\citep{gunter2021compromised}.

\begin{figure}[t!]
  \centering
    \subfloat[$(-0.3,-0.7)$]
  {\includegraphics[width=0.21\textwidth]{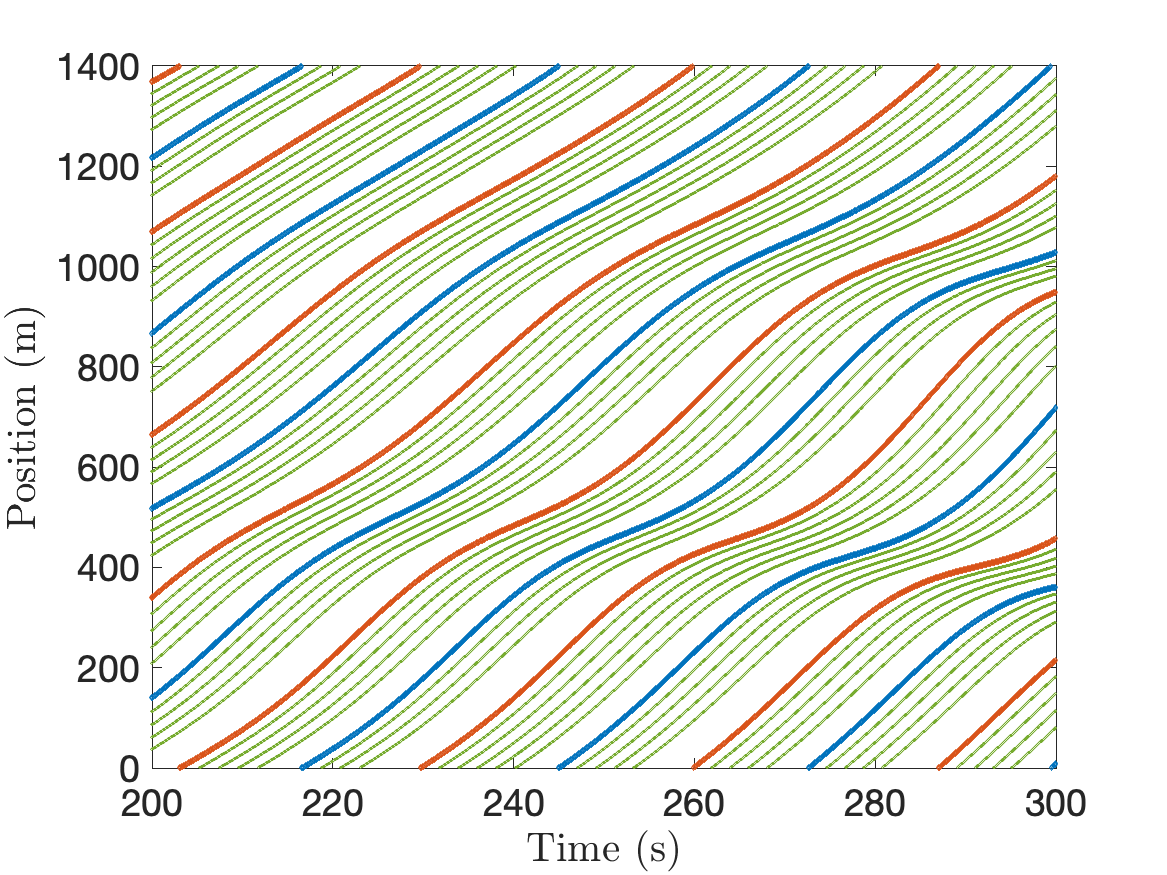}\label{trajectory_type_II_01}}
  \hspace*{-0.8em}%
    \subfloat[$(-0.3,-0.3)$]{\includegraphics[width=0.21\textwidth]{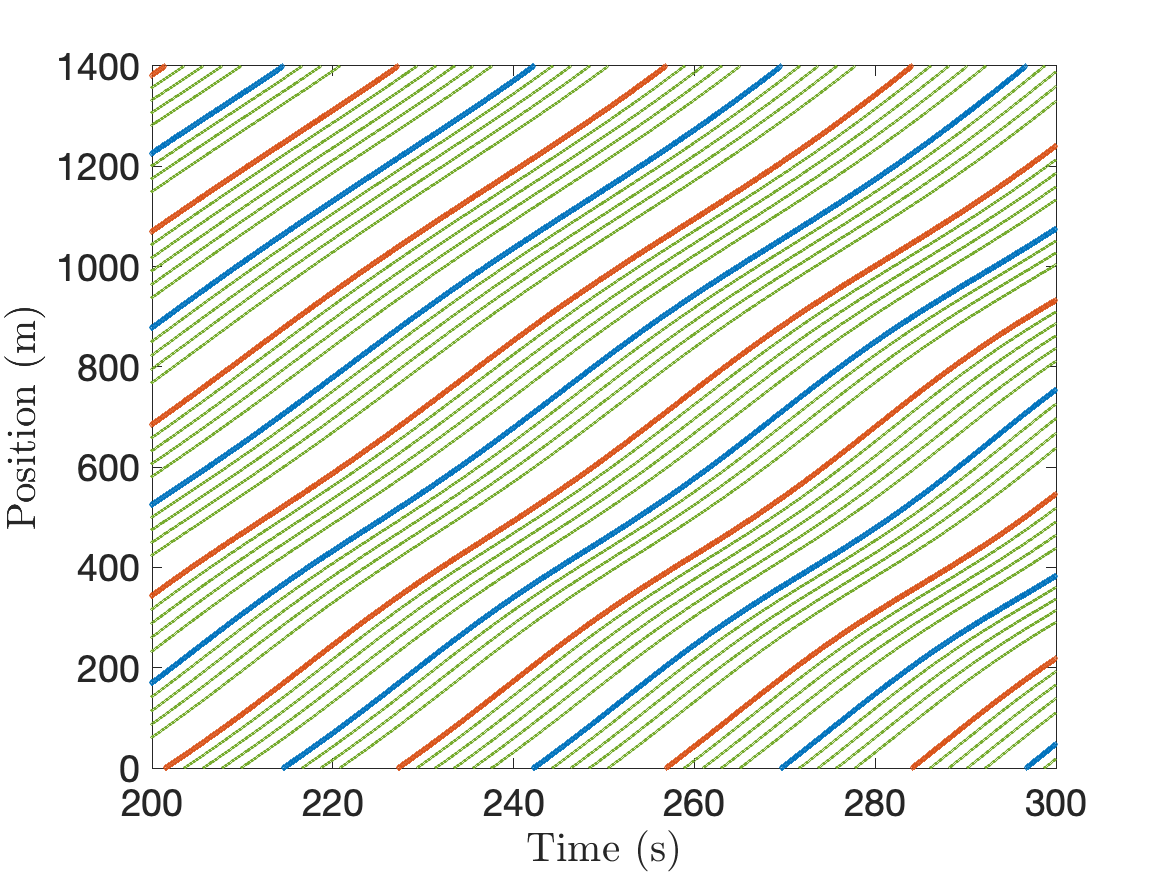}\label{trajectory_type_II_02}}
  \hspace*{-0.8em}%
    \subfloat[$(-0.3,0.0)$]{\includegraphics[width=0.21\textwidth]{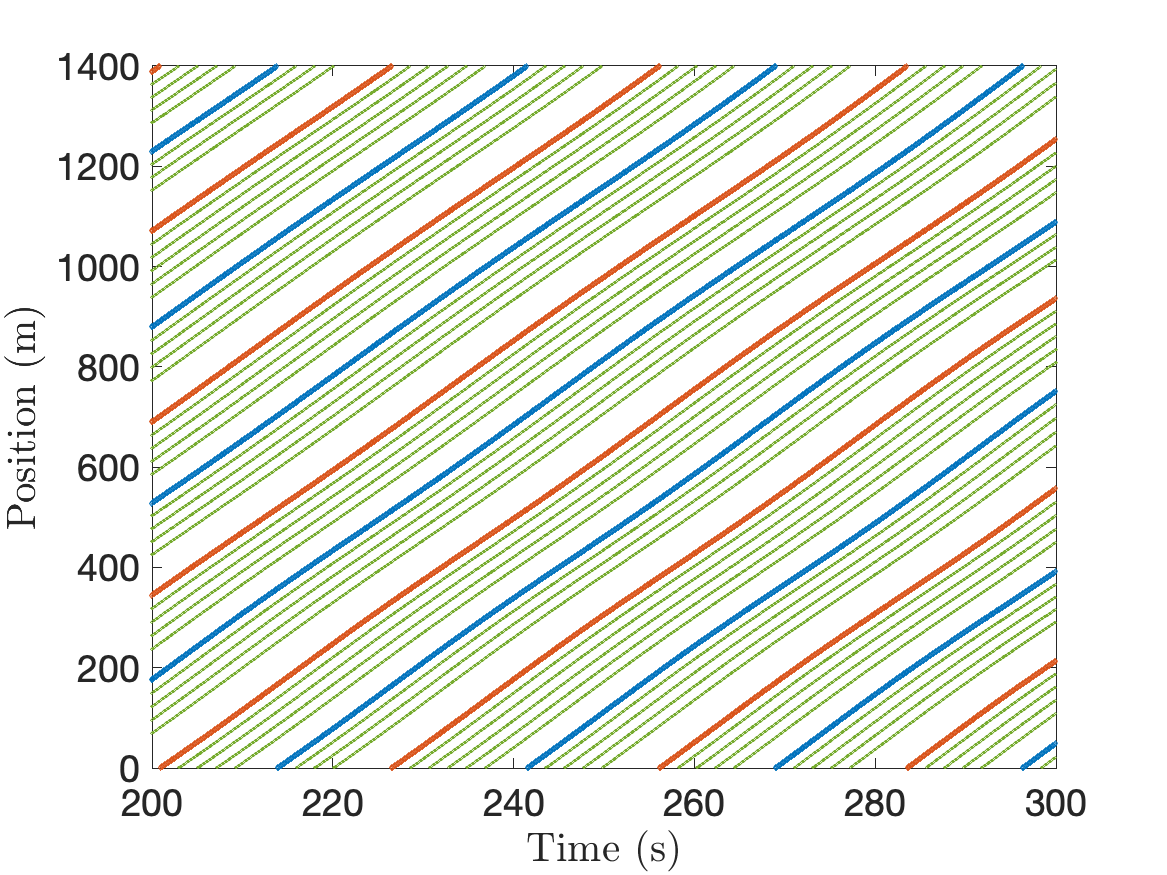}\label{trajectory_type_II_03}}
  \hspace*{-0.8em}%
    \subfloat[$(-0.3,0.4)$]{\includegraphics[width=0.21\textwidth]{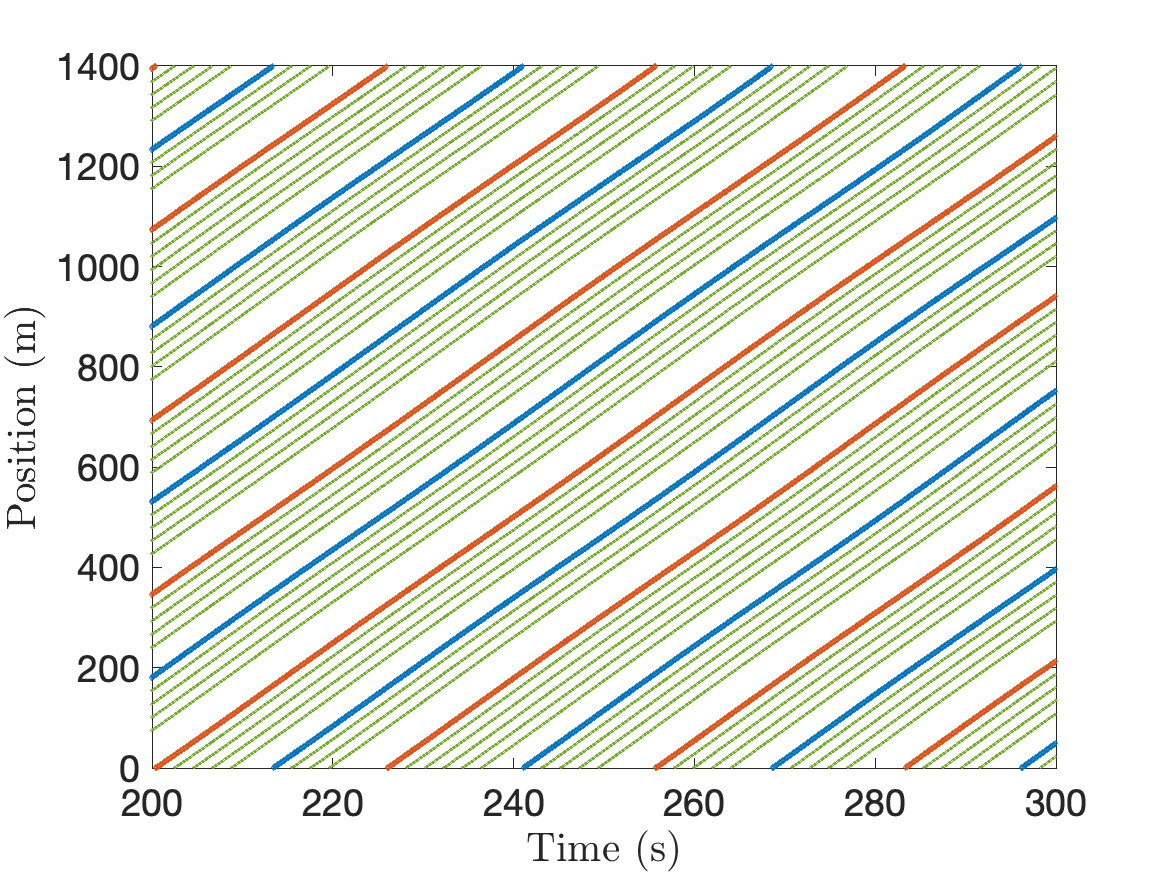}\label{trajectory_type_II_04}}
  \hspace*{-0.8em}%
    \subfloat[$(-0.3,0.8)$]{\includegraphics[width=0.21\textwidth]{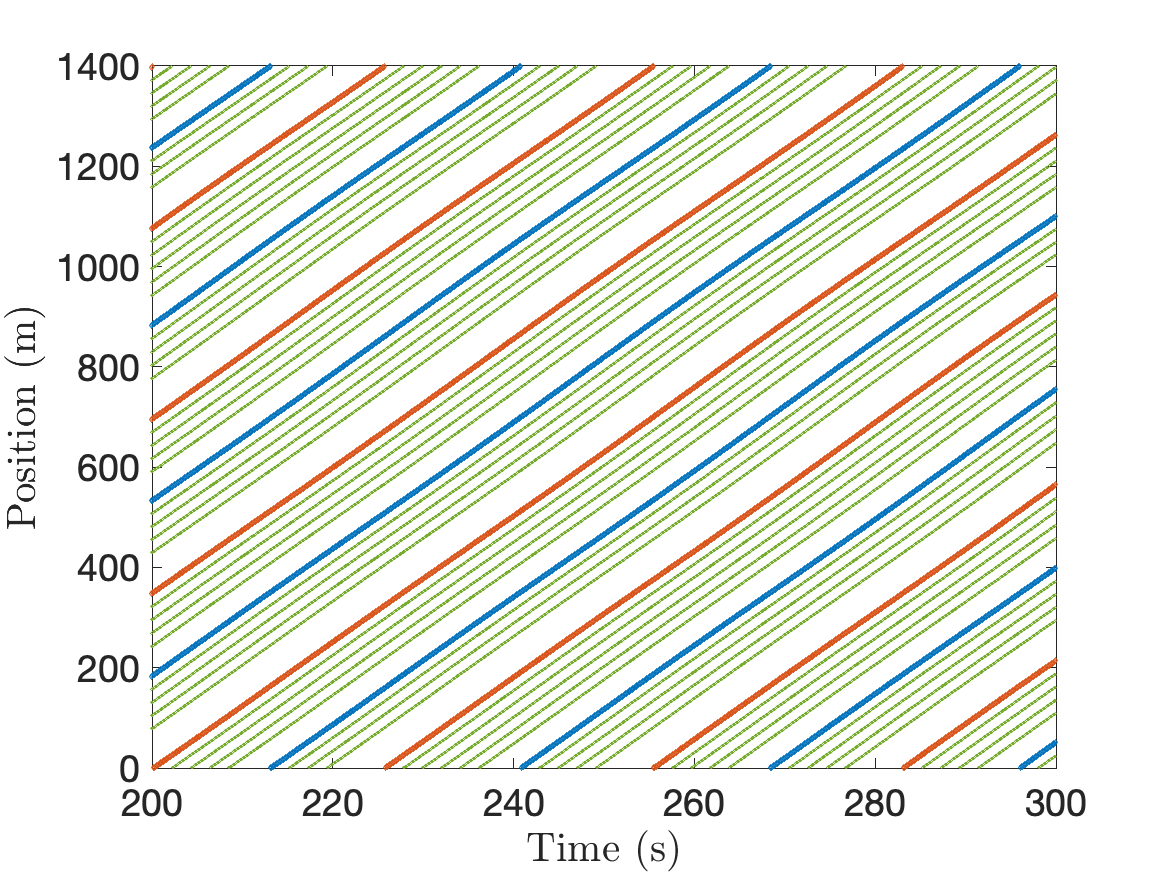}\label{trajectory_type_II_05}}%
  \vfil%
  \vspace*{-1.0em}%
    \subfloat[$(0.0,-0.7)$]{\includegraphics[width=0.21\textwidth]{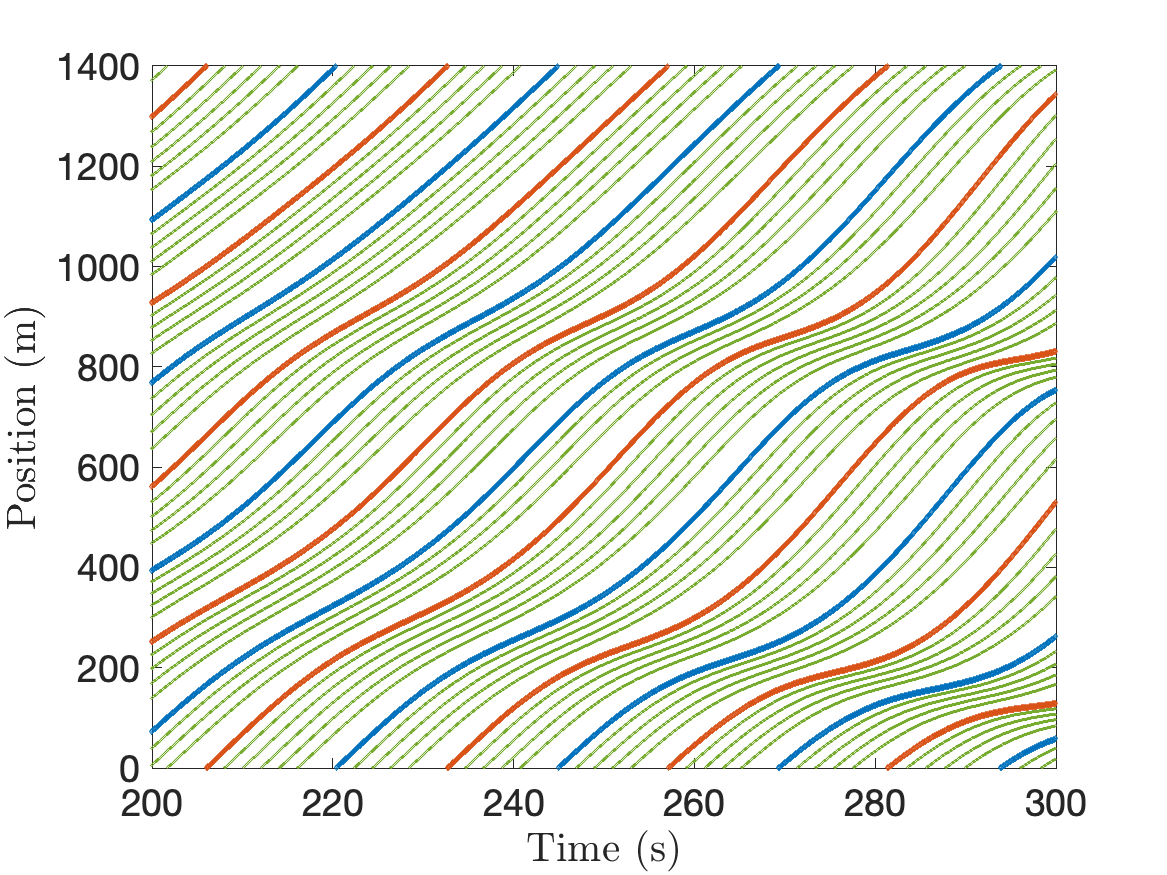}\label{trajectory_type_II_06}}
  \hspace*{-0.8em}%
    \subfloat[$(0.0,-0.3)$]{\includegraphics[width=0.21\textwidth]{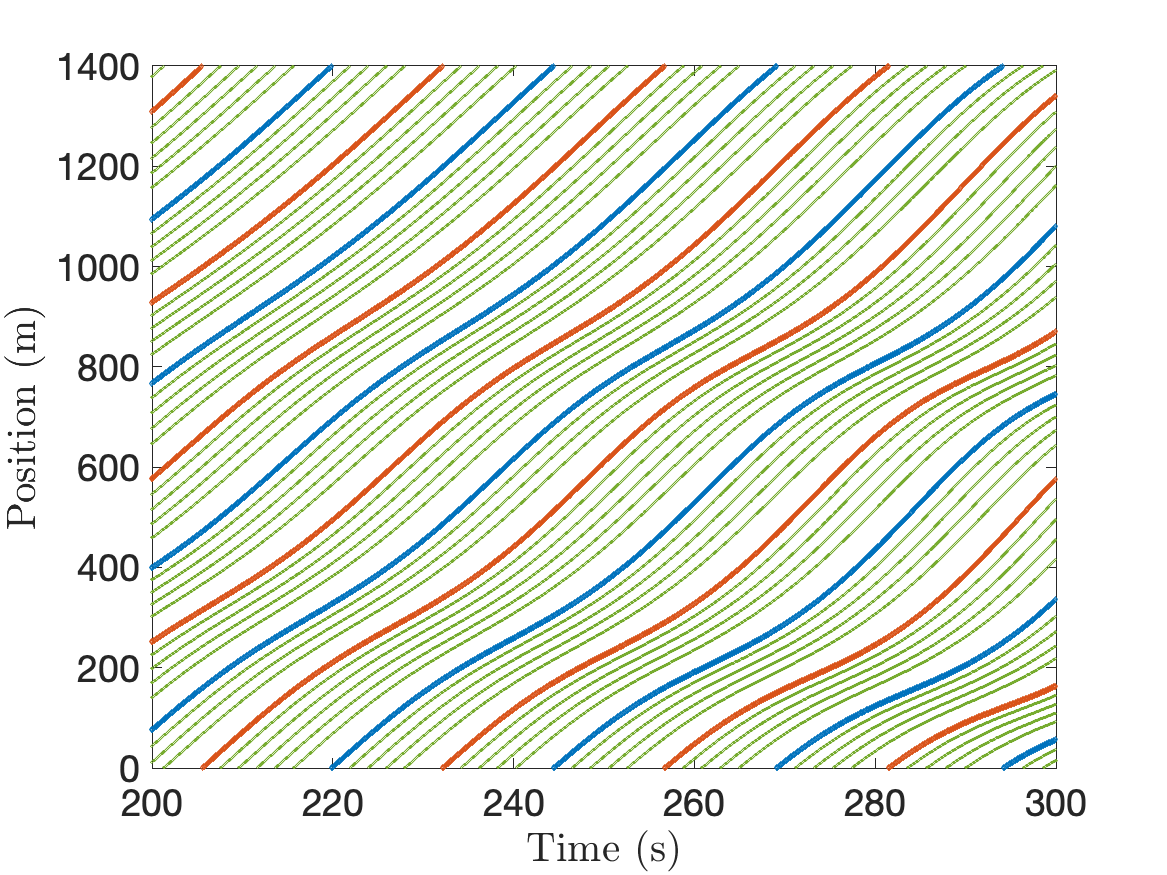}\label{trajectory_type_II_07}}
  \hspace*{-0.8em}%
    \subfloat[$(0.0,0.0)$]{\includegraphics[width=0.21\textwidth]{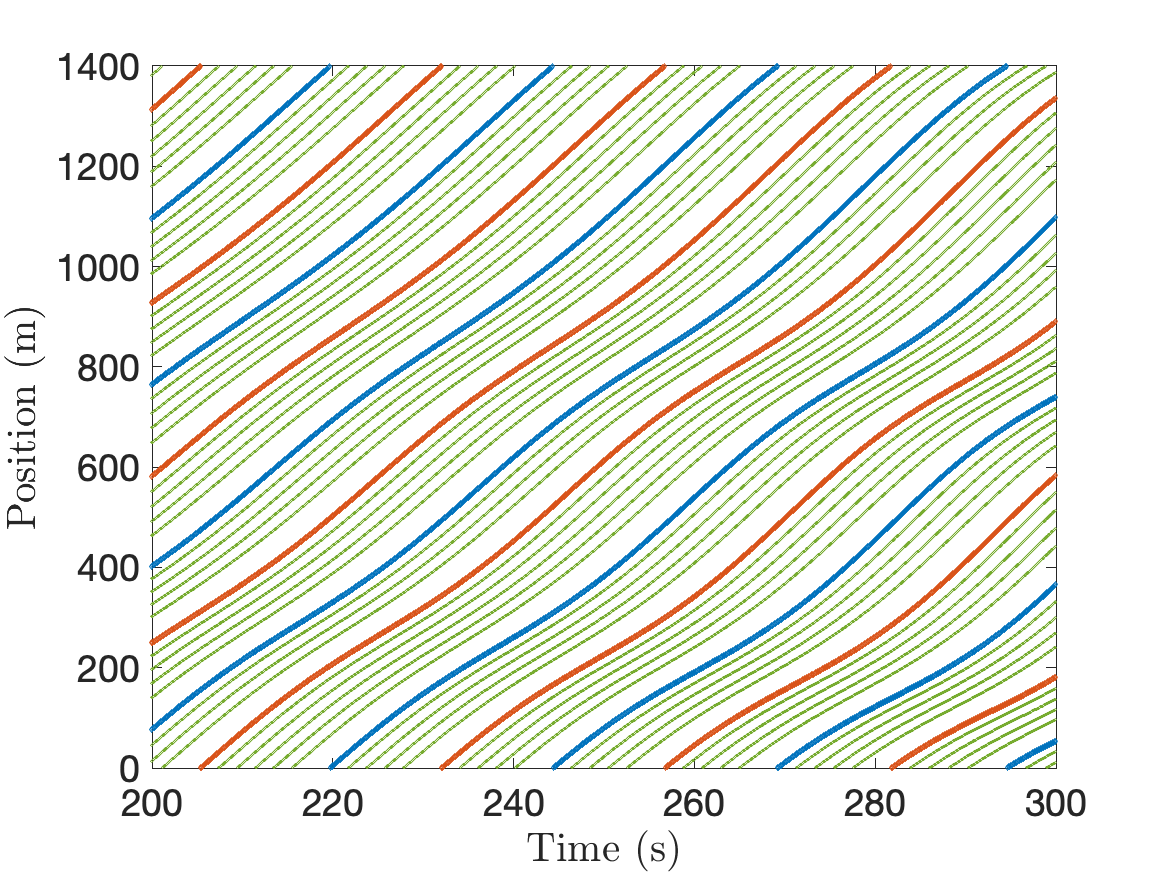}\label{trajectory_type_II_08}}
  \hspace*{-0.8em}%
    \subfloat[$(0.0,0.4)$]{\includegraphics[width=0.21\textwidth]{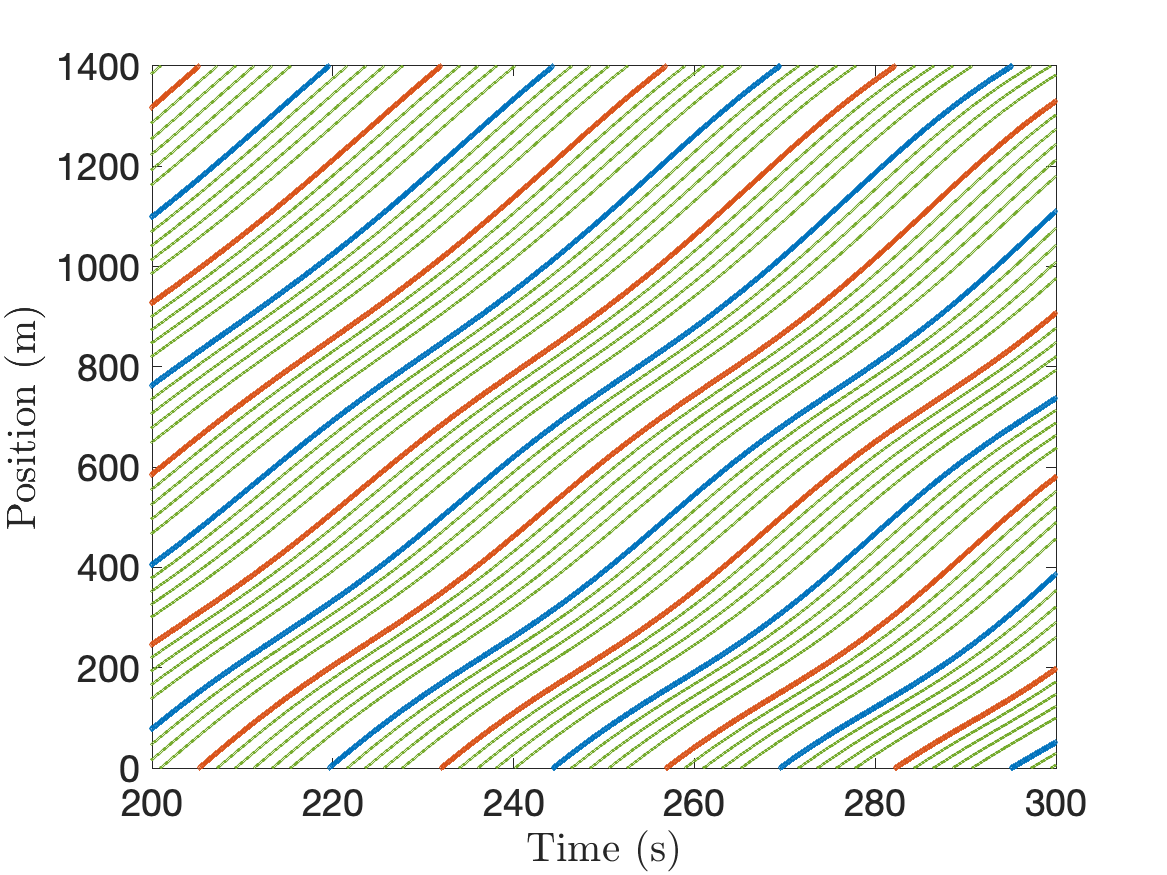}\label{trajectory_type_II_09}}
  \hspace*{-0.8em}%
    \subfloat[$(0.0,0.8)$]{\includegraphics[width=0.21\textwidth]{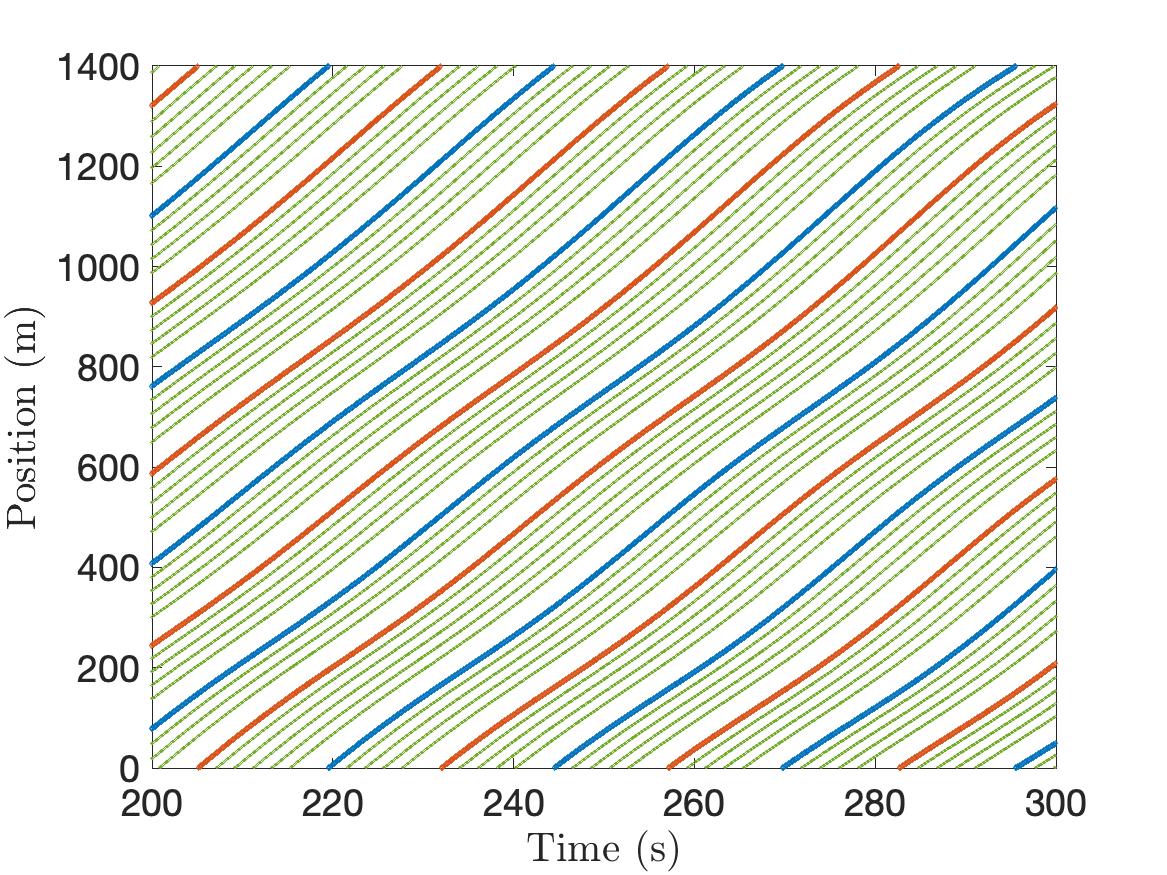}\label{trajectory_type_II_10}}%
  \vfil%
  \vspace*{-1.0em}%
    \subfloat[$(0.4,-0.7)$]{\includegraphics[width=0.21\textwidth]{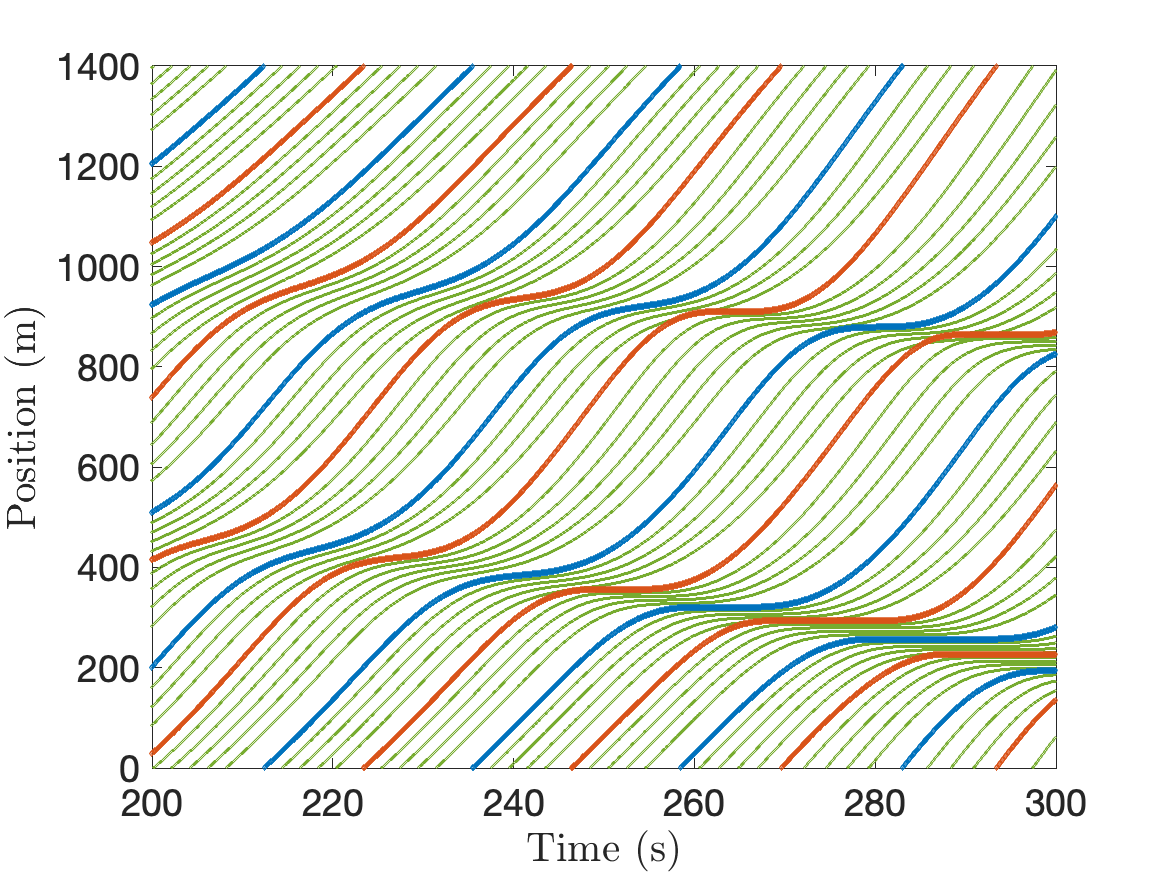}\label{trajectory_type_II_11}}
  \hspace*{-0.8em}%
    \subfloat[$(0.4,-0.3)$]{\includegraphics[width=0.21\textwidth]{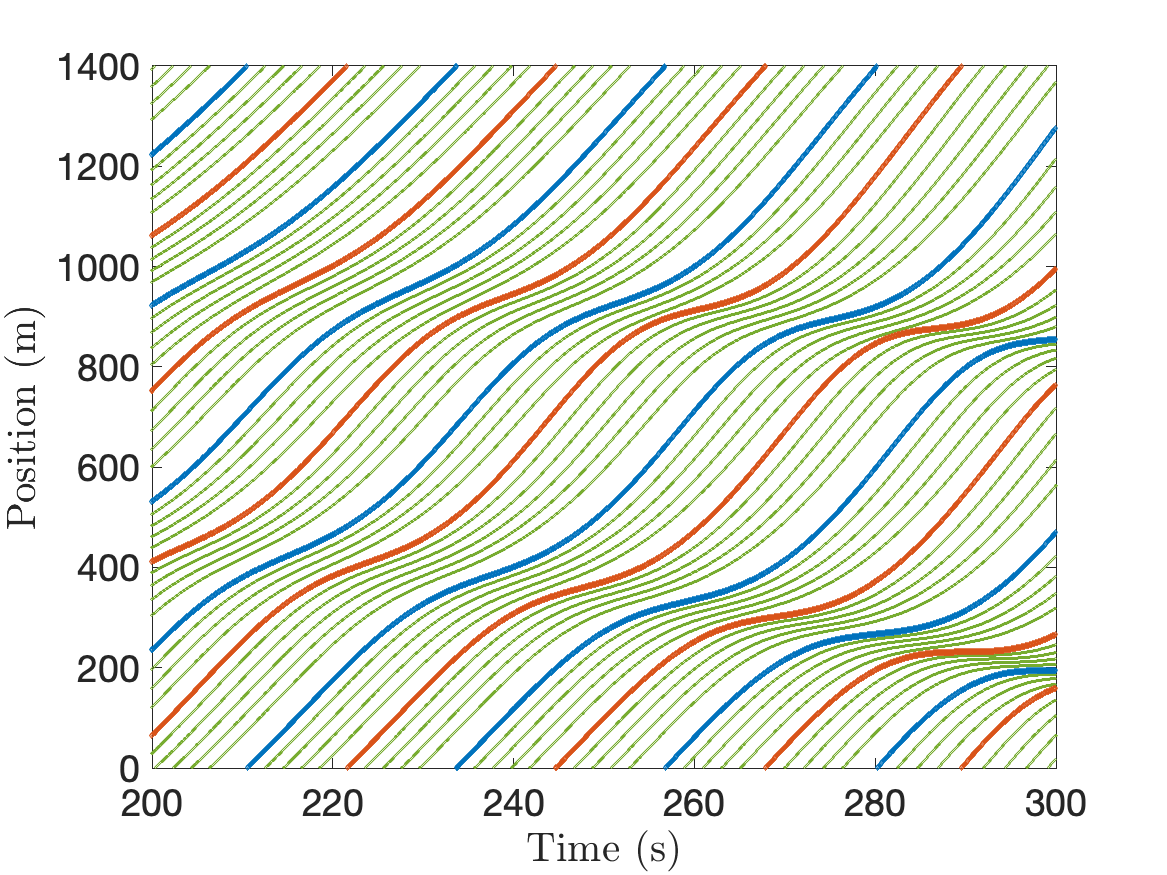}\label{trajectory_type_II_12}}
  \hspace*{-0.8em}%
    \subfloat[$(0.4,0.0)$]{\includegraphics[width=0.21\textwidth]{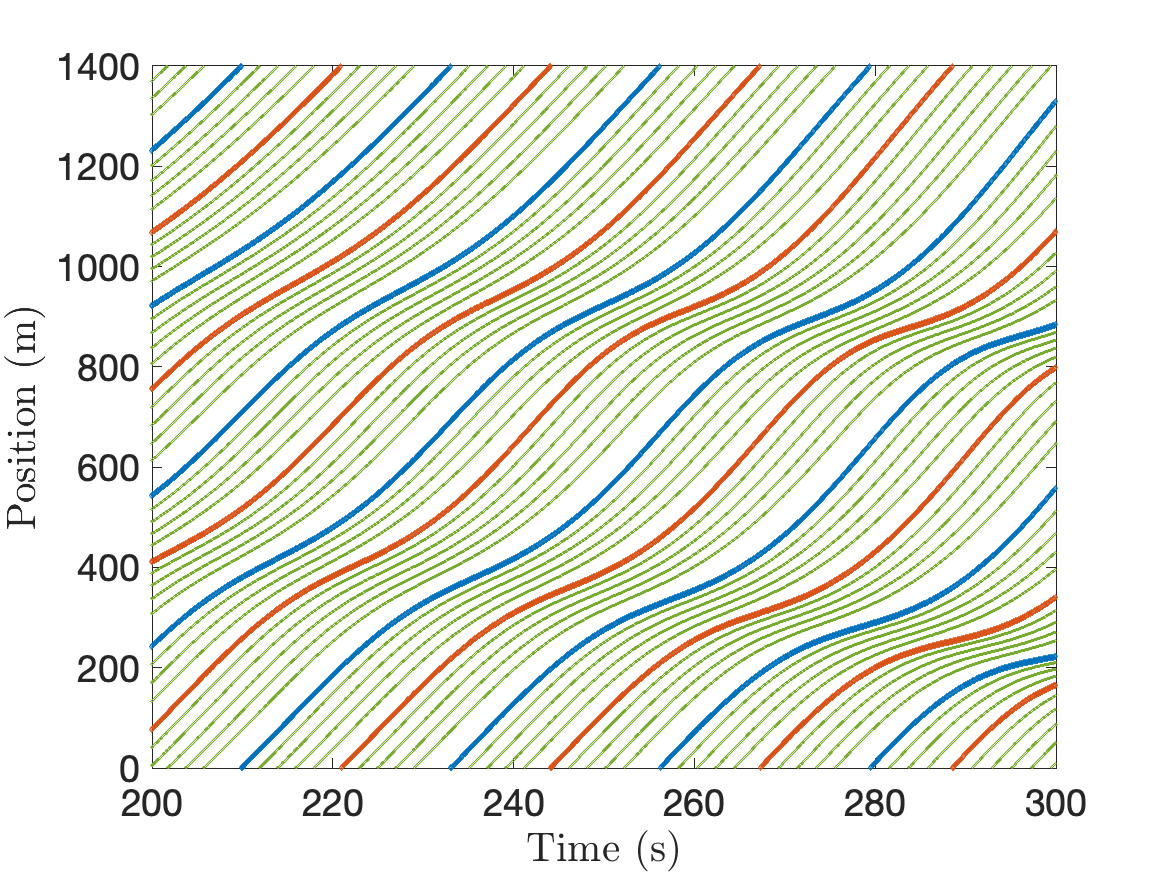}\label{trajectory_type_II_13}}
  \hspace*{-0.8em}%
    \subfloat[$(0.4,0.4)$]{\includegraphics[width=0.21\textwidth]{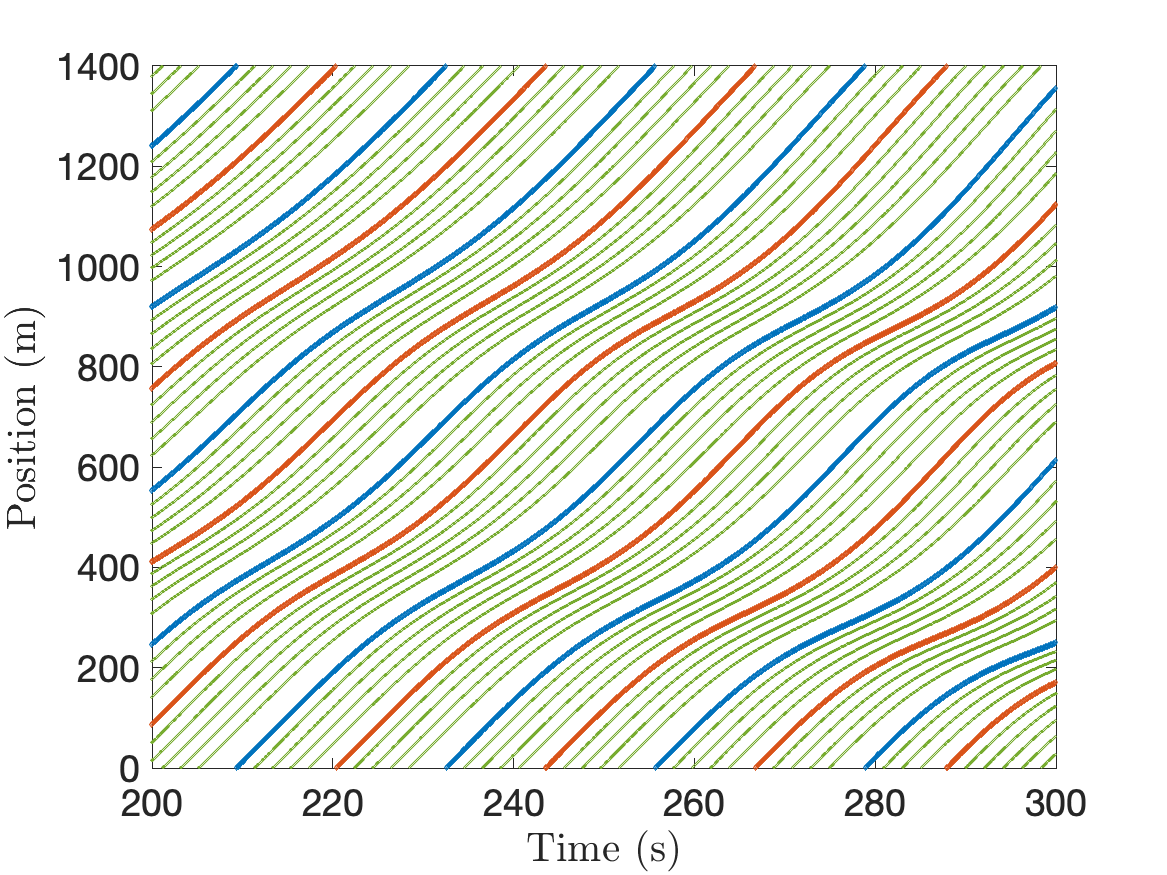}\label{trajectory_type_II_14}}
  \hspace*{-0.8em}%
    \subfloat[$(0.4,0.8)$]{\includegraphics[width=0.21\textwidth]{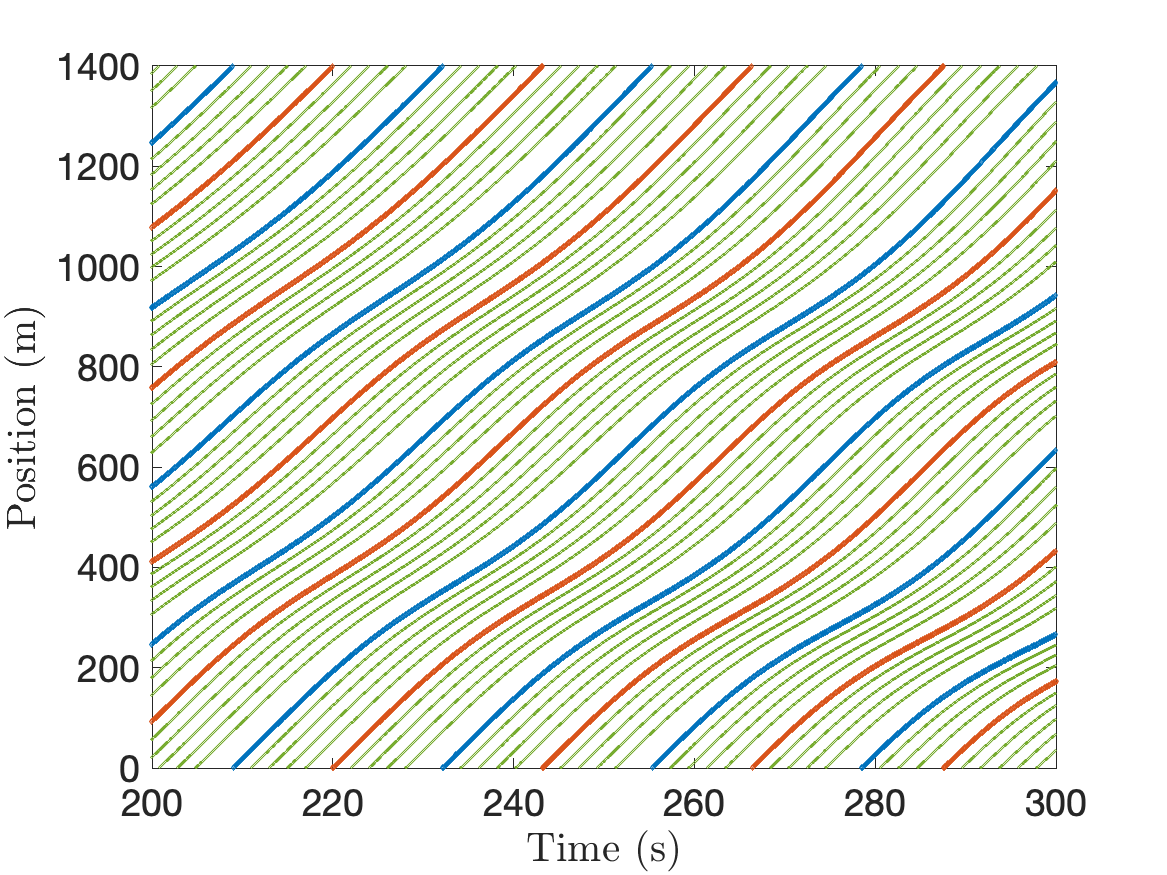}\label{trajectory_type_II_15}}%
  \vfil%
  \vspace*{-1.0em}%
    \subfloat[$(0.8,-0.7)$]{\includegraphics[width=0.21\textwidth]{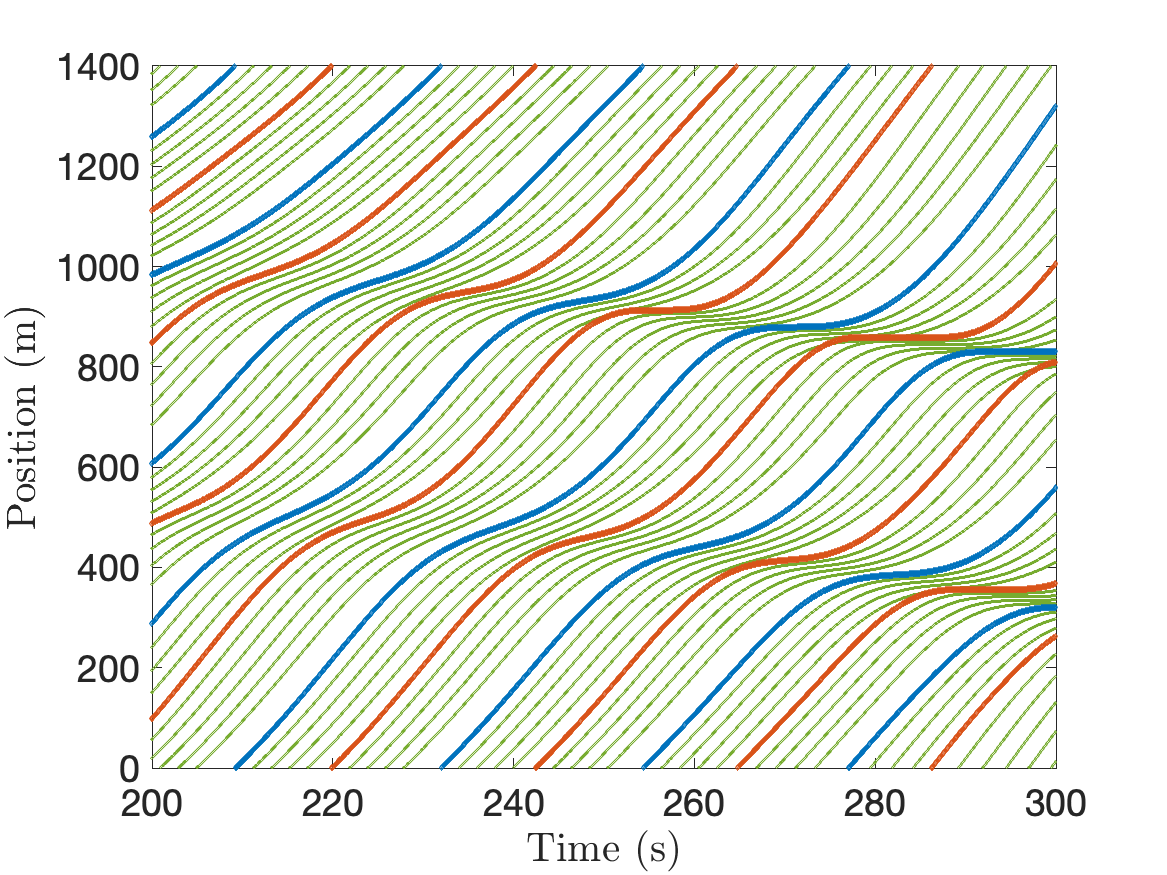}\label{trajectory_type_II_16}}
  \hspace*{-0.8em}
    \subfloat[$(0.8,-0.3)$]{\includegraphics[width=0.21\textwidth]{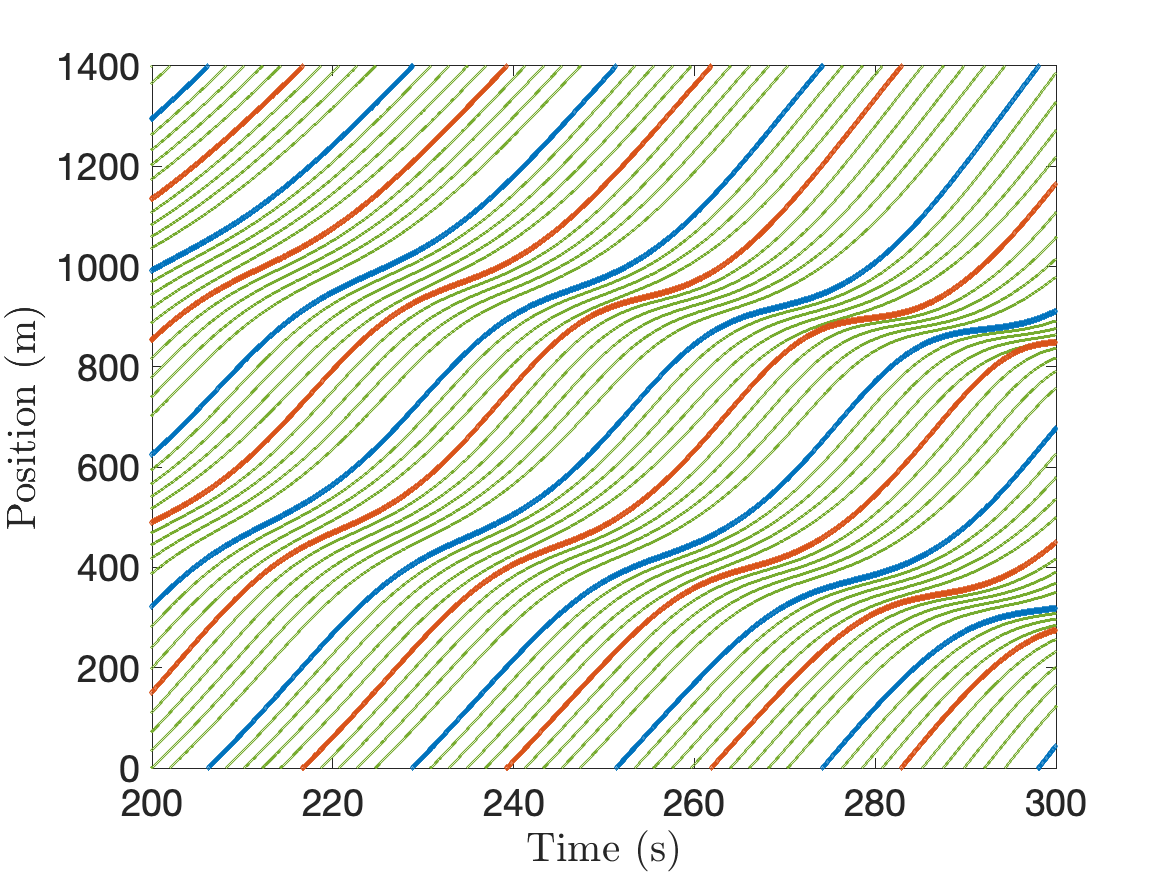}\label{trajectory_type_II_17}}
  \hspace*{-0.8em}%
    \subfloat[$(0.8,0.0)$]{\includegraphics[width=0.21\textwidth]{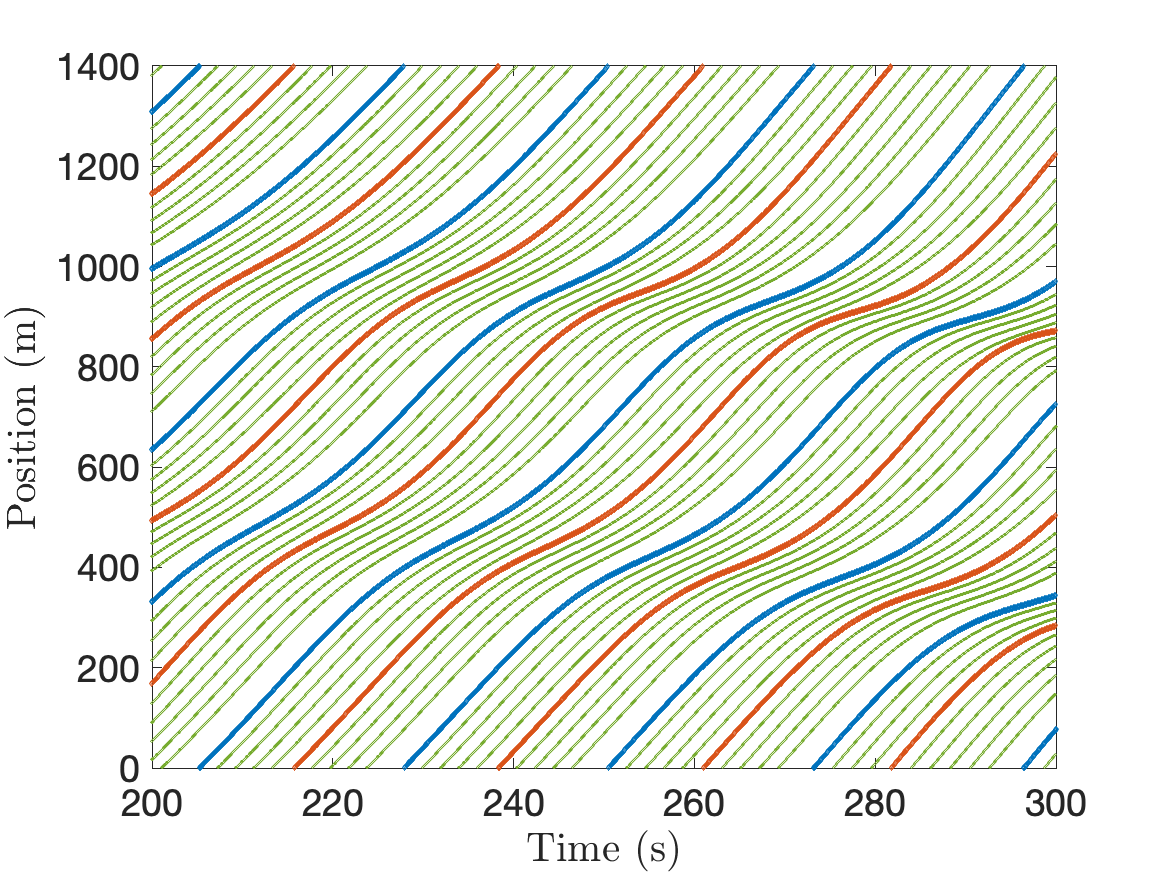}\label{trajectory_type_II_18}}
  \hspace*{-0.8em}%
    \subfloat[$(0.8,0.4)$]{\includegraphics[width=0.21\textwidth]{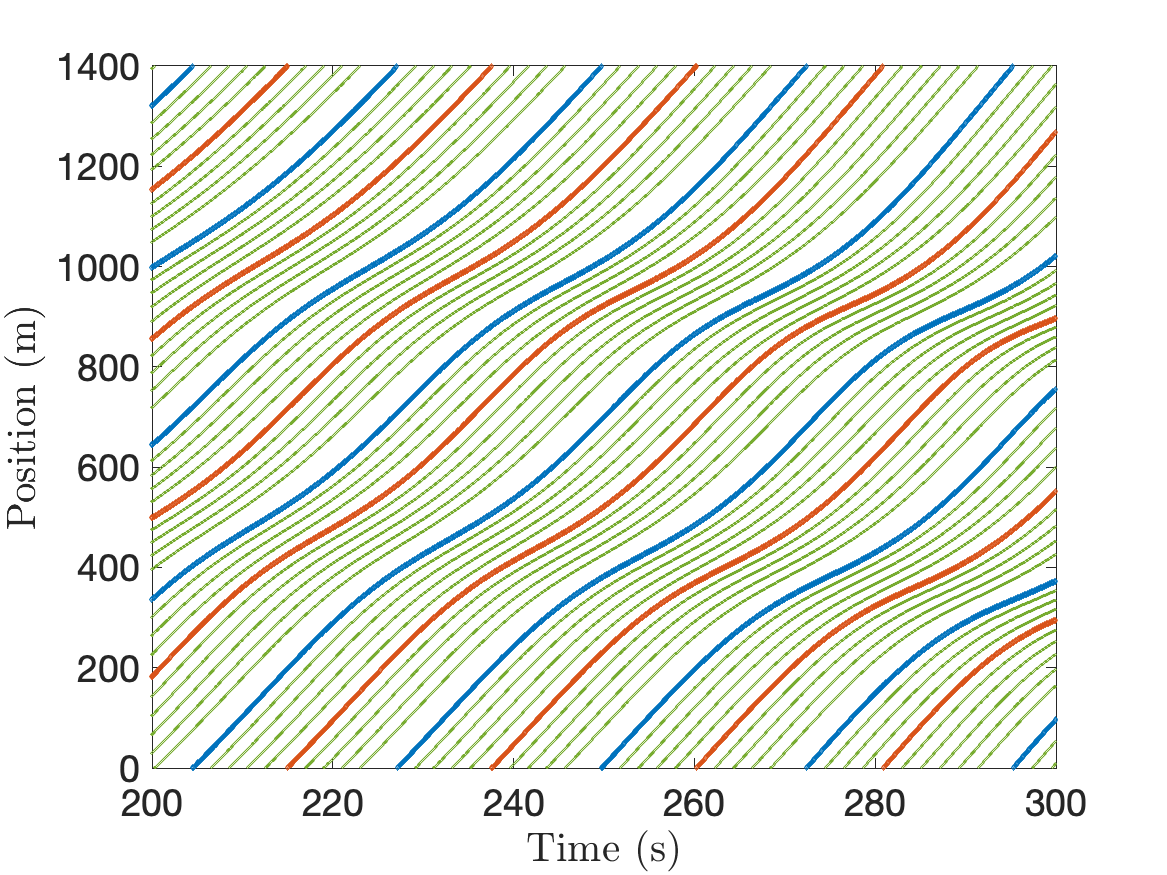}\label{trajectory_type_II_19}}
  \hspace*{-0.8em}%
    \subfloat[$(0.8,0.8)$]{\includegraphics[width=0.21\textwidth]{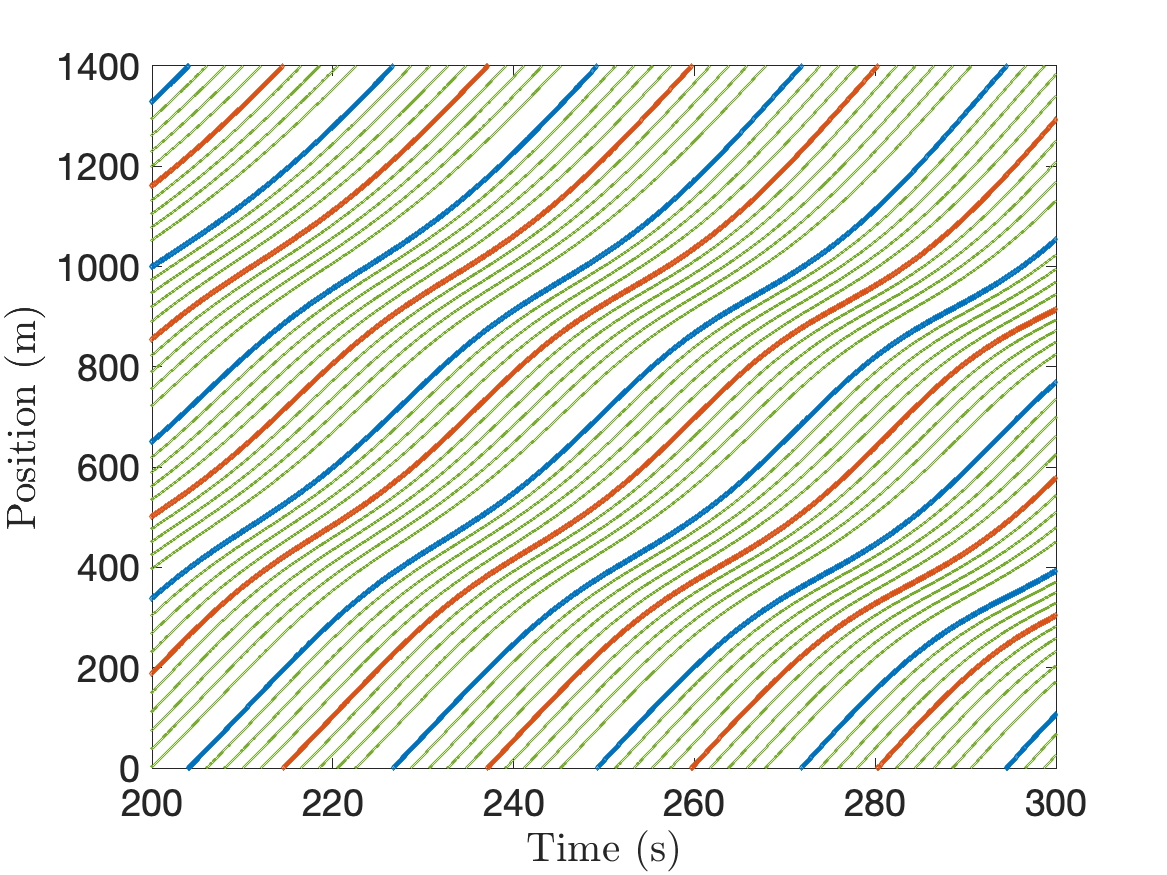}\label{trajectory_type_II_20}}
  \caption{Trajectory at different pairs of $(\delta_{1,i},\delta_{2,i})$. Green, blue, and red correspond to HDV, normal ACC vehicle, and attacked ACC vehicle, respectively.}\label{trajectory_type_II}%
\end{figure}

\subsection{Type~II attacks}\label{section5.3}

Following the conditions given by equation~\eqref{eq4.14}, we use the values of $\delta_{1,i} \in [-0.3, 0.8]$ and $\delta_{2,i} \in [-0.7, 0.8]$ for Type~II candidate attacks.  

\subsubsection{Trajectories and safety impacts}

Following the same simulation setting as in Section~\ref{section5.2}, vehicle trajectories for Type~II attacks are shown in Fig.~\ref{trajectory_type_II}, with values of the pair $(\delta_{1,i}, \delta_{2,i})$ given above for candidate attacks. Fig.~\ref{trajectory_type_II_08} shows vehicle trajectories in the absence of attack, serving as a baseline for comparison.

When solely increasing (or decreasing) the value of $\delta_{1,i}$, more (or less) oscillations are observed in the traffic flow. This is because a positive $\delta_{1,i}$ fools the ACC vehicle to believe in a greater spacing gap to its preceding vehicle, thereby resulting in a larger acceleration. Consequently, the attacked vehicle drives in a more aggressive manner with a smaller spacing gap to its lead vehicle; and vice versa for a negative $\delta_{1,i}$. Similarly, when increasing $\delta_{2,i}$ alone, fewer oscillations are observed due to improvement of the $\hat{\lambda}_2$ value as indicated by equation~\eqref{eq4.6}. These results are consistent with the findings that a larger positive $\delta_{1,i}$ and a smaller negative $\delta_{2,i}$ tend to result in a greater $\hat{\lambda}_2$, as shown in Fig.~\ref{lambda2-hat}.

\subsubsection{Traffic efficiency (fundamental diagram)}

\begin{figure}[t!]
  \centering
  \includegraphics[width=\textwidth]{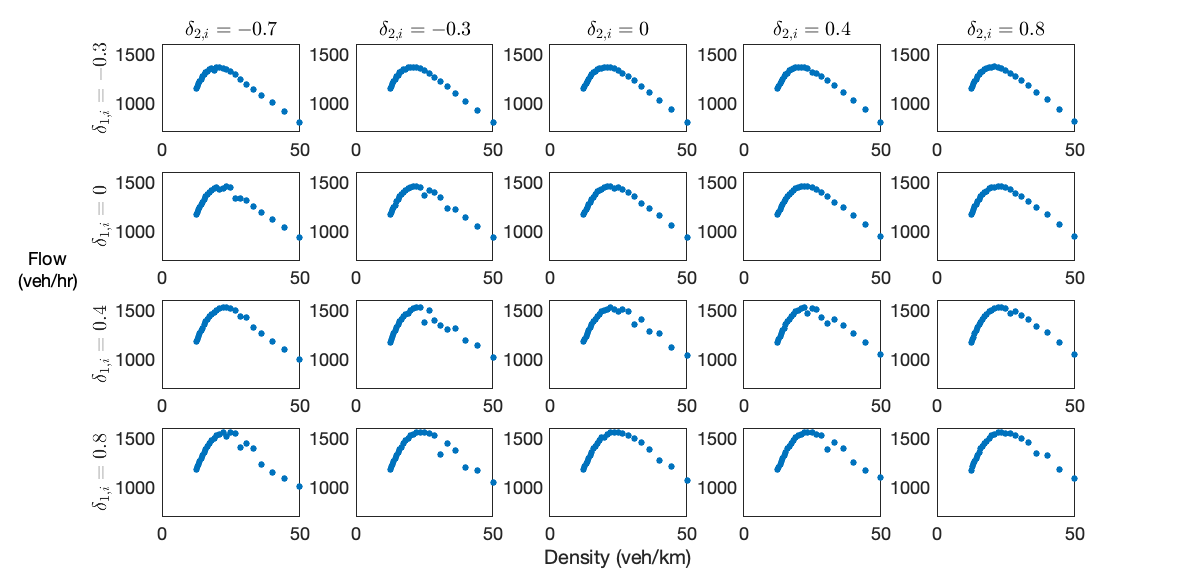}
  \vspace*{-2.0em}%
  \caption{Fundamental diagrams for different values of the pair $(\delta_{1,i}, \delta_{2,i})$ considering Type~II attacks.}
  \label{FD_type_II}
\end{figure}

To assess traffic efficiency in the presence of Type~II attacks, we present the corresponding FD in Fig.~\ref{FD_type_II}. The results show that traffic capacity increases with the increase of $\delta_{1,i}$ for a certain $\delta_{2,i}$. This is due to the fact that a larger $\delta_{1,i}$ would lead to greater vehicle acceleration, resulting in smaller spacing gap and thereby increasing the throughput. However, when the value of $\delta_{1,i}$ is fixed and $\delta_{2,i}$ increases from -0.7 to 0.8, traffic capacity remains almost unchanged, indicating that attacks on the relative speed $\Delta v$ may have less impact on traffic efficiency compared to those on the spacing $s$.

Compared with the non-attack scenario, increasing the value of $\delta_{1,i}$ may increase traffic capacity. However, more traffic oscillations would arise due to the resulting aggressive driving behavior of compromised ACC vehicles. In contrast, it is noted that increasing the value of $\delta_{2,i}$ may not significantly impact traffic efficiency, but could result in less traffic oscillations.

\begin{remark}\label{remark5.1}
It is worth noting that $\lambda_2 > 0$ in the absence of attack for all the scenarios considered. Hence, leveraging the conditions of~\eqref{eq3.16} and~\eqref{eq4.14} allows one to simulate candidate attacks that destabilize ACC dynamics as well as those that degrade ACC stability.
\end{remark}

\section{Concluding Remarks}\label{section6}

In this study, we have developed a general framework for modeling and synthesizing two types of candidate attacks on ACC vehicles with analytical characterization of their mathematical properties, covering direct attacks on ACC control commands and false data injection attacks on ACC sensor measurement. Under the analytical framework developed, we derive explicit conditions to characterize the malicious nature of potential attacks. Those conditions provide an effective method for synthesizing an array of strategic attacks on ACC vehicles, with a higher degree of realism in modeling their adverse effects compared to many prior studies. Extensive numerical simulations have been conducted to illustrate the mechanism of candidate attacks, offering useful insights into understanding the vulnerability of future transportation systems. For example, positive candidate attacks (in terms of $\delta_i$ and $\delta_{1,i}$ for Type~I and II attacks, respectively) likely lead to larger vehicle acceleration resulting in more aggressive driving behavior of the attacked vehicle, thereby compromising traffic safety. In contrast, negative candidate attacks (in terms of $\delta_i$ and $\delta_{2,i}$ for Type~I and II attacks, respectively) tend to significantly degrade traffic efficiency without causing much safety concerns. 

Following many prior studies, we have used the widely adopted OVRV model for describing the car-following dynamics of ACC vehicles. While it has been shown to fit well to both simulated and real trajectories of vehicles equipped with ACC capabilities, it is not exactly clear how highly automated vehicles are to be operated in the future. Hence, one interesting direction for further research is to use reinforcement learning techniques to synthesize model-free cyberattacks. That is, malicious actors could alter the driving behavior of a compromised vehicle in a rational manner (for remaining stealthy) based on the real traffic data accessed, to degrade the performance of transportation systems. While it is impossible to obtain the exact form of all potential attacks, this study provides an effective method for synthesizing candidate attacks with a good degree of realism, thereby offers useful insights into understanding their widespread impacts on transportation cybersecurity. This is expected to motivate a line of research on the development of effective attack detection and mitigation strategies in future studies, such as anomaly detection of traffic flow and advanced control design of intelligent vehicles.

\bibliographystyle{unsrtnat}
\bibliography{mybibfile}

\begin{thebibliography}{45}
\providecommand{\natexlab}[1]{#1}
\providecommand{\url}[1]{\texttt{#1}}
\expandafter\ifx\csname urlstyle\endcsname\relax
  \providecommand{\doi}[1]{doi: #1}\else
  \providecommand{\doi}{doi: \begingroup \urlstyle{rm}\Url}\fi

\bibitem[Sun et~al.(2022)Sun, Wang, Shao, Sun, and Levin]{sun2022energy}
Wenbo Sun, Shian Wang, Yunli Shao, Zongxuan Sun, and Michael~W Levin.
\newblock Energy and mobility impacts of connected autonomous vehicles with
  co-optimization of speed and powertrain on mixed vehicle platoons.
\newblock \emph{Transportation Research Part C: Emerging Technologies},
  142:\penalty0 103764, 2022.

\bibitem[Wang et~al.(2022{\natexlab{a}})Wang, Stern, and
  Levin]{wang2022optimal}
Shian Wang, Raphael Stern, and Michael~W Levin.
\newblock Optimal control of autonomous vehicles for traffic smoothing.
\newblock \emph{IEEE Transactions on Intelligent Transportation Systems},
  23\penalty0 (4):\penalty0 3842--3852, 2022{\natexlab{a}}.

\bibitem[Wang et~al.(2021)Wang, Levin, and Caverly]{wang2021optimal}
Shian Wang, Michael~W Levin, and Ryan~James Caverly.
\newblock Optimal parking management of connected autonomous vehicles: A
  control-theoretic approach.
\newblock \emph{Transportation Research Part C: Emerging Technologies},
  124:\penalty0 102924, 2021.

\bibitem[Parkinson et~al.(2017)Parkinson, Ward, Wilson, and
  Miller]{parkinson2017cyber}
Simon Parkinson, Paul Ward, Kyle Wilson, and Jonathan Miller.
\newblock Cyber threats facing autonomous and connected vehicles: Future
  challenges.
\newblock \emph{IEEE Transactions on Intelligent Transportation Systems},
  18\penalty0 (11):\penalty0 2898--2915, 2017.

\bibitem[Dong et~al.(2020)Dong, Wang, Ni, Liu, and Chen]{dong2020impact}
Changyin Dong, Hao Wang, Daiheng Ni, Yongfei Liu, and Quan Chen.
\newblock Impact evaluation of cyber-attacks on traffic flow of connected and
  automated vehicles.
\newblock \emph{IEEE Access}, 8:\penalty0 86824--86835, 2020.

\bibitem[Petit and Shladover(2014)]{petit2014potential}
Jonathan Petit and Steven Shladover.
\newblock Potential cyberattacks on automated vehicles.
\newblock \emph{IEEE Transactions on Intelligent Transportation Systems},
  16\penalty0 (2):\penalty0 546--556, 2014.

\bibitem[Alipour-Fanid et~al.(2020)Alipour-Fanid, Dabaghchian, and
  Zeng]{alipour2020impact}
Amir Alipour-Fanid, Monireh Dabaghchian, and Kai Zeng.
\newblock Impact of jamming attacks on vehicular cooperative adaptive cruise
  control systems.
\newblock \emph{IEEE Transactions on Vehicular Technology}, 69\penalty0
  (11):\penalty0 12679--12693, 2020.

\bibitem[Li et~al.(2018)Li, Tu, Fan, Dong, and Wang]{li2018influence}
Ye~Li, Yu~Tu, Qi~Fan, Changyin Dong, and Wei Wang.
\newblock Influence of cyber-attacks on longitudinal safety of connected and
  automated vehicles.
\newblock \emph{Accident Analysis \& Prevention}, 121:\penalty0 148--156, 2018.

\bibitem[Li et~al.(2023)Li, Rosenblad, Wang, Shang, and Stern]{li2023exploring}
Tianyi Li, Benjamin Rosenblad, Shian Wang, Mingfeng Shang, and Raphael Stern.
\newblock Exploring energy impacts of cyberattacks on adaptive cruise control
  vehicles.
\newblock In \emph{2023 IEEE Intelligent Vehicles Symposium (IV)}, pages 1--6.
  IEEE, 2023.

\bibitem[Li et~al.(2022)Li, Shang, Wang, Filippelli, and
  Stern]{li2022detecting}
Tianyi Li, Mingfeng Shang, Shian Wang, Matthew Filippelli, and Raphael Stern.
\newblock Detecting stealthy cyberattacks on automated vehicles via generative
  adversarial networks.
\newblock In \emph{IEEE 25th International Conference on Intelligent
  Transportation Systems}, pages 3632--3637. IEEE, 2022.

\bibitem[Zhou et~al.(2022)Zhou, Wang, and Peeta]{zhou2022robust}
Anye Zhou, Jian Wang, and Srinivas Peeta.
\newblock Robust control strategy for platoon of connected and autonomous
  vehicles considering falsified information injected through communication
  links.
\newblock \emph{Journal of Intelligent Transportation Systems}, pages 1--17,
  2022.

\bibitem[Wang et~al.(2023{\natexlab{a}})Wang, Levin, and Stern]{wang2023minmax}
Shian Wang, Michael~W Levin, and Raphael Stern.
\newblock Optimal feedback control law for automated vehicles in the presence
  of cyberattacks: A min--max approach.
\newblock \emph{Transportation Research Part C: Emerging Technologies},
  153:\penalty0 104204, 2023{\natexlab{a}}.

\bibitem[Wang et~al.(2020)Wang, Wu, and He]{wang2020modeling}
Pengcheng Wang, Xinkai Wu, and Xiaozheng He.
\newblock Modeling and analyzing cyberattack effects on connected automated
  vehicular platoons.
\newblock \emph{Transportation Research Part C: Emerging Technologies},
  115:\penalty0 102625, 2020.

\bibitem[Wang et~al.(2023{\natexlab{b}})Wang, Zhang, Masoud, and
  Liu]{wang2023anomaly}
Yiyang Wang, Ruixuan Zhang, Neda Masoud, and Henry~X Liu.
\newblock Anomaly detection and string stability analysis in connected
  automated vehicular platoons.
\newblock \emph{Transportation Research Part C: Emerging Technologies},
  151:\penalty0 104114, 2023{\natexlab{b}}.

\bibitem[Zeadally et~al.(2012)Zeadally, Hunt, Chen, Irwin, and
  Hassan]{zeadally2012vehicular}
Sherali Zeadally, Ray Hunt, Yuh-Shyan Chen, Angela Irwin, and Aamir Hassan.
\newblock Vehicular ad hoc networks {(VANETS)}: status, results, and
  challenges.
\newblock \emph{Telecommunication Systems}, 50\penalty0 (4):\penalty0 217--241,
  2012.

\bibitem[Wang and Chen(2021)]{wang2021resilient}
Fengchen Wang and Yan Chen.
\newblock Resilient flocking control for connected and automated vehicles with
  cyber-attack threats.
\newblock \emph{ASME Letters in Dynamic Systems and Control}, 1\penalty0
  (3):\penalty0 1--6, 2021.

\bibitem[Wang(2023, \text{to appear})]{wang2023novel}
Shian Wang.
\newblock A novel framework for modeling and synthesizing stealthy cyberattacks
  on driver-assist enabled vehicles.
\newblock In \emph{2023 IEEE Intelligent Vehicles Symposium}. IEEE, 2023,
  \text{to appear}.

\bibitem[Khraisat et~al.(2019)Khraisat, Gondal, Vamplew, and
  Kamruzzaman]{khraisat2019survey}
Ansam Khraisat, Iqbal Gondal, Peter Vamplew, and Joarder Kamruzzaman.
\newblock Survey of intrusion detection systems: techniques, datasets and
  challenges.
\newblock \emph{Cybersecurity}, 2\penalty0 (1):\penalty0 1--22, 2019.

\bibitem[Wang et~al.(2022{\natexlab{b}})Wang, Li, and
  Levin]{wang2022optimalTRC}
Shian Wang, Zhexian Li, and Michael~W Levin.
\newblock Optimal policy for integrating autonomous vehicles into the auto
  market.
\newblock \emph{Transportation Research Part C: Emerging Technologies},
  143:\penalty0 103821, 2022{\natexlab{b}}.

\bibitem[Wilson and Ward(2011)]{wilson2011car}
R~Eddie Wilson and Jonathan~A Ward.
\newblock Car-following models: Fifty years of linear stability analysis--a
  mathematical perspective.
\newblock \emph{Transportation Planning and Technology}, 34\penalty0
  (1):\penalty0 3--18, 2011.

\bibitem[Ye et~al.(2023)Ye, Sun, and Sun]{ye2022car}
Yingjun Ye, Jie Sun, and Jian Sun.
\newblock Car-following characteristics of commercially available adaptive
  cruise control systems and comparison with human drivers.
\newblock \emph{Transportation Research Record}, 2677\penalty0 (2):\penalty0
  1401--1414, 2023.

\bibitem[Talebpour and Mahmassani(2016)]{talebpour2016influence}
Alireza Talebpour and Hani~S Mahmassani.
\newblock Influence of connected and autonomous vehicles on traffic flow
  stability and throughput.
\newblock \emph{Transportation Research Part C: Emerging Technologies},
  71:\penalty0 143--163, 2016.

\bibitem[Wang et~al.(2023{\natexlab{c}})Wang, Shang, Levin, and
  Stern]{wang2023general}
Shian Wang, Mingfeng Shang, Michael~W Levin, and Raphael Stern.
\newblock A general approach to smoothing nonlinear mixed traffic via control
  of autonomous vehicles.
\newblock \emph{Transportation Research Part C: Emerging Technologies},
  146:\penalty0 103967, 2023{\natexlab{c}}.

\bibitem[Rajamani(2011)]{rajamani2011vehicle}
Rajesh Rajamani.
\newblock \emph{Vehicle dynamics and control}.
\newblock Springer Science \& Business Media, 2011.

\bibitem[Wang(2022)]{wang2022planning}
Shian Wang.
\newblock \emph{Planning, operation, and management of automated transportation
  systems: A control-theoretic approach}.
\newblock PhD thesis, University of Minnesota, 2022.

\bibitem[van Wyk et~al.(2019)van Wyk, Wang, Khojandi, and Masoud]{van2019real}
Franco van Wyk, Yiyang Wang, Anahita Khojandi, and Neda Masoud.
\newblock Real-time sensor anomaly detection and identification in automated
  vehicles.
\newblock \emph{IEEE Transactions on Intelligent Transportation Systems},
  21\penalty0 (3):\penalty0 1264--1276, 2019.

\bibitem[Boem et~al.(2017)Boem, Gallo, Ferrari-Trecate, and
  Parisini]{boem2017distributed}
Francesca Boem, Alexander~J Gallo, Giancarlo Ferrari-Trecate, and Thomas
  Parisini.
\newblock A distributed attack detection method for multi-agent systems
  governed by consensus-based control.
\newblock In \emph{2017 IEEE 56th Annual Conference on Decision and Control
  (CDC)}, pages 5961--5966. IEEE, 2017.

\bibitem[Biron et~al.(2018)Biron, Dey, and Pisu]{biron2018real}
Zoleikha~Abdollahi Biron, Satadru Dey, and Pierluigi Pisu.
\newblock Real-time detection and estimation of denial of service attack in
  connected vehicle systems.
\newblock \emph{IEEE Transactions on Intelligent Transportation Systems},
  19\penalty0 (12):\penalty0 3893--3902, 2018.

\bibitem[Treiber et~al.(2000)Treiber, Hennecke, and
  Helbing]{treiber2000congested}
Martin Treiber, Ansgar Hennecke, and Dirk Helbing.
\newblock Congested traffic states in empirical observations and microscopic
  simulations.
\newblock \emph{Physical Review E}, 62\penalty0 (2):\penalty0 1805--1824, 2000.

\bibitem[Sarker et~al.(2019)Sarker, Shen, Rahman, Chowdhury, Dey, Li, Wang, and
  Narman]{sarker2019review}
Ankur Sarker, Haiying Shen, Mizanur Rahman, Mashrur Chowdhury, Kakan Dey,
  Fangjian Li, Yue Wang, and Husnu~S Narman.
\newblock A review of sensing and communication, human factors, and controller
  aspects for information-aware connected and automated vehicles.
\newblock \emph{IEEE Transactions on Intelligent Transportation Systems},
  21\penalty0 (1):\penalty0 7--29, 2019.

\bibitem[Pourabdollah et~al.(2017)Pourabdollah, Bj{\"a}rkvik, F{\"u}rer,
  Lindenberg, and Burgdorf]{pourabdollah2017calibration}
Mitra Pourabdollah, Eric Bj{\"a}rkvik, Florian F{\"u}rer, Bj{\"o}rn Lindenberg,
  and Klaas Burgdorf.
\newblock Calibration and evaluation of car following models using real-world
  driving data.
\newblock In \emph{2017 IEEE 20th International Conference on Intelligent
  Transportation Systems}, pages 1--6. IEEE, 2017.

\bibitem[He and Wang(2023)]{he2023calibrating}
Linjia He and Xuesong Wang.
\newblock Calibrating car-following models on urban streets using naturalistic
  driving data.
\newblock \emph{Journal of Transportation Engineering, Part A: Systems},
  149\penalty0 (4):\penalty0 1--16, 2023.

\bibitem[Milan{\'e}s et~al.(2013)Milan{\'e}s, Shladover, Spring,
  et~al.]{milanes2013cooperative}
Vicente Milan{\'e}s, Steven~E Shladover, John Spring, et~al.
\newblock Cooperative adaptive cruise control in real traffic situations.
\newblock \emph{IEEE Transactions on Intelligent Transportation Systems},
  15\penalty0 (1):\penalty0 296--305, 2013.

\bibitem[Ioannou and Chien(1993)]{ioannou1993autonomous}
Petros~A Ioannou and Cheng-Chih Chien.
\newblock Autonomous intelligent cruise control.
\newblock \emph{IEEE Transactions on Vehicular Technology}, 42\penalty0
  (4):\penalty0 657--672, 1993.

\bibitem[Montanino and Punzo(2021)]{montanino2021string}
Marcello Montanino and Vincenzo Punzo.
\newblock On string stability of a mixed and heterogeneous traffic flow: A
  unifying modelling framework.
\newblock \emph{Transportation Research Part B: Methodological}, 144:\penalty0
  133--154, 2021.

\bibitem[Sugiyama et~al.(2008)Sugiyama, Fukui, Kikuchi,
  et~al.]{sugiyama2008traffic}
Yuki Sugiyama, Minoru Fukui, Macoto Kikuchi, et~al.
\newblock Traffic jams without bottlenecks-experimental evidence for the
  physical mechanism of the formation of a jam.
\newblock \emph{New Journal of Physics}, 10\penalty0 (3):\penalty0 033001,
  2008.

\bibitem[Wu et~al.(2019)Wu, Stern, Cui, et~al.]{wu2019tracking}
Fangyu Wu, Raphael~E Stern, Shumo Cui, et~al.
\newblock Tracking vehicle trajectories and fuel rates in phantom traffic jams:
  Methodology and data.
\newblock \emph{Transportation Research Part C: Emerging Technologies},
  99:\penalty0 82--109, 2019.

\bibitem[Treiber and Kesting(2013)]{treiber2013traffic}
Martin Treiber and Arne Kesting.
\newblock \emph{Traffic flow dynamics}.
\newblock Springer, 2013.

\bibitem[Gunter et~al.(2019)Gunter, Stern, and Work]{gunter2019modeling}
George Gunter, Raphael Stern, and Daniel~B Work.
\newblock Modeling adaptive cruise control vehicles from experimental data:
  model comparison.
\newblock In \emph{2019 IEEE Intelligent Transportation Systems Conference},
  pages 3049--3054. IEEE, 2019.

\bibitem[Hayward(1972)]{hayward1972near}
John~C Hayward.
\newblock Near miss determination through use of a scale of danger.
\newblock 1972.

\bibitem[Li et~al.(2020)Li, Wu, Lee, Yang, and Shi]{li2020analysis}
Ye~Li, Dan Wu, Jaeyoung Lee, Min Yang, and Yuntao Shi.
\newblock Analysis of the transition condition of rear-end collisions using
  time-to-collision index and vehicle trajectory data.
\newblock \emph{Accident Analysis \& Prevention}, 144:\penalty0 105676, 2020.

\bibitem[Mahmud et~al.(2019)Mahmud, Ferreira, Hoque, and
  Tavassoli]{mahmud2019micro}
SM~Sohel Mahmud, Luis Ferreira, Md~Shamsul Hoque, and Ahmad Tavassoli.
\newblock Micro-simulation modelling for traffic safety: A review and potential
  application to heterogeneous traffic environment.
\newblock \emph{IATSS Research}, 43\penalty0 (1):\penalty0 27--36, 2019.

\bibitem[Greenshields et~al.(1934)Greenshields, Thompson, Dickinson, and
  Swinton]{greenshields1934photographic}
Bruce~Douglas Greenshields, J~T Thompson, H~C Dickinson, and R~S Swinton.
\newblock The photographic method of studying traffic behavior.
\newblock In \emph{Highway Research Board Proceedings}, volume~13, 1934.

\bibitem[Coifman(2015)]{coifman2015empirical}
Benjamin Coifman.
\newblock Empirical flow-density and speed-spacing relationships: Evidence of
  vehicle length dependency.
\newblock \emph{Transportation Research Part B: Methodological}, 78:\penalty0
  54--65, 2015.

\bibitem[Gunter et~al.(2021)Gunter, Li, Hojjati, et~al.]{gunter2021compromised}
George Gunter, Huichen Li, Avesta Hojjati, et~al.
\newblock Compromised {ACC} vehicles can degrade current mixed-autonomy traffic
  performance while remaining stealthy against detection.
\newblock \emph{arXiv preprint arXiv:2112.11986}, 2021.

\end{thebibliography}

\end{document}